\documentclass[ number, 5p, twocolumn,sort&compress,times,fleqn]{elsarticle}

\usepackage{amsmath}
\usepackage{amssymb}
\usepackage{graphicx}
\usepackage{subfig}
\usepackage{bm}
\usepackage{siunitx}
\usepackage[colorlinks = true,
            linkcolor = blue,
            urlcolor  = blue,
            citecolor = blue,
            anchorcolor = blue]{hyperref}

\usepackage{./csml}
\usepackage{soul}

\usepackage{mathtools}
\graphicspath{{./figs/}}

\begin{document}

\begin{frontmatter}

\title{Subdivision surfaces with isogeometric analysis adapted refinement weights} 

\author[cam]{Qiaoling Zhang}
\ead{zq217@cam.ac.uk}

\author[cam,num]{Malcolm Sabin}
\ead{malcolm.sabin@btinternet.com}

\author[cam]{Fehmi Cirak\corref{cor1}}
\ead{f.cirak@eng.cam.ac.uk}

\cortext[cor1]{Corresponding author}

\address[cam]{Department of Engineering,  University of Cambridge, 
Trumpington Street, Cambridge CB2 1PZ, UK}
\address[num]{Numerical Geometry Ltd, 19 John Amner Close, Ely, Cambridgeshire CB6 1DT, UK}


\begin{abstract}
Subdivision surfaces provide an elegant isogeometric analysis framework for geometric design and analysis of partial differential equations defined on surfaces. They are already a standard in high-end computer animation and graphics and are becoming available in a number of geometric modelling systems for engineering design. The subdivision refinement rules  are usually adapted from knot insertion rules for splines. The quadrilateral Catmull-Clark scheme considered in this work is equivalent to cubic B-splines away from extraordinary, or irregular, vertices with other than four adjacent elements.  Around extraordinary vertices the surface consists of a nested sequence of smooth spline patches which join $C^1$ continuously at the point itself. As known from geometric design literature, the subdivision weights can be optimised so that the surface quality is improved by minimising short-wavelength surface oscillations around extraordinary vertices. We use the related techniques to determine weights that minimise finite element discretisation errors as measured in the thin-shell energy norm. The optimisation problem is formulated over a characteristic domain and the errors in approximating cup- and saddle-like quadratic shapes obtained from eigenanalysis of the subdivision matrix are minimised. In finite element analysis the optimised subdivision weights for either cup- or saddle-like shapes are chosen depending on the shape of the solution field around an extraordinary vertex. As our computations confirm, the optimised subdivision weights yield a reduction of $50\%$ and more in discretisation errors in the energy and $L_2$ norms. Although, as to be expected, the convergence rates are the same as for the classical Catmull-Clark weights, the convergence constants are improved.

\end{abstract}

\newpage


\end{frontmatter}

%
\section{Introduction \label{sec: introduction}}
Isogeometric analysis aims to provide a seamless engineering design-analysis workflow by using a single common representation for geometric modelling and analysis.  This is usually achieved by  representing geometry and discretising analysis models with the same kind of basis functions~\cite{Hughes:2005aa}. The prevailing feature-based CAD modelling systems rely on trimmed NURBS and boundary representations (B-Reps). The resulting non-watertight geometries consisting of several trimmed patches pose unique challenges to finite element analysis. As a generalisation of splines, subdivision surfaces can provide watertight representations for geometries with arbitrary topology. After their early success in computer animation and graphics they are now supported in many CAD systems, including Catia, PTC Creo and Autodesk Fusion 360. Before the advent of isogeometric analysis, it had already been realised that subdivision surfaces provide also ideal basis functions for finite element analysis, in particular, of thin-shells~\cite{Cirak:2000aa,reif2001curvature, Cirak:2001aa, Cirak:2002aa}, see also more recent work~\cite{Cirak:2011aa,bandara2018isogeometric}. 

Subdivision schemes for generating smooth surfaces were first described in the late 1970s as an extension of low degree B-splines to control meshes with non-tensor-product connectivity~\cite{Doo:1978aa, Catmull:1978aa}. In subdivision a geometry is described with a control mesh and a limiting process of repeated refinement. For parts of the mesh containing only regular vertices, with each adjacent to four quadrilateral faces, the refinement rules are adapted from knot insertion rules for B-splines. For the remaining parts with extraordinary vertices the refinement rules are chosen such that they yield in the limit a smooth surface. Subdivision refinement is a linear mapping of coordinates of the coarse control mesh  to the coordinates of the refined mesh with a subdivision matrix. Hence, the local limit surface properties can be inferred from the eigenstructure of the subdivision matrix after a discrete Fourier transform~\cite{Doo:1978aa,Reif:1995aa}. The $C^1$ continuity of the surface and its curvature behaviour at the extraordinary vertex depend on eigenvalues and the ordering, i.e. Fourier indices, of the corresponding eigenvectors. In turn, both depend on the coefficients of the subdivision matrix that encodes the specific refinement rules applied. 

As known, around extraordinary vertices short-wavelength surface oscillations, i.e. ripples, may occur irrespective of $C^1$ continuity and boundedness of curvature~\cite{peters2004shape, karvciauskas2004shape}. There have been many attempts to improve the fairness of subdivision surfaces, that is, to minimise curvature variations, by carefully tuning the refinement rules, earlier works include~\cite{halstead1993efficient, kobbelt1996variational}. More recently, in Augsd\"orfer et al.~\cite{augsdorfer2006tuning} the refinement rules for Catmull-Clark and other quadrilateral schemes have been optimised such that the variation of the Gaussian curvature is minimised while ensuring bounded curvatures. Different from the direct search method used in~\cite{augsdorfer2006tuning}, the refinement rules can also be obtained from a nonlinear constrained optimisation problem. Barthe et al.~\cite{barthe2004subdivision} apply such a procedure to triangular Loop and $\sqrt{3}$-subdivision schemes with a multi-objective cost function comprised of terms penalising divergence of curvatures and aiming local quadratic precision. In Ginkel et al.~\cite{ginkel2008local} a fairness increasing cost function containing the third derivatives of the surface in combination with $C^1$ continuity and bounded curvature constraints is optimised. 

In the present paper, we optimise the subdivision refinement rules so that their approximation properties are improved when used in finite element analysis of thin-shells. Thin-shells are prevalent in many engineering applications, most prominently in aerospace, automotive and structural engineering, and are equivalent to thin-plates when their unstressed geometry is planar~\cite{Ciarlet:2005aa}. The thin-shell energy functional, and weak form, depend on the second order derivatives of the stressed surface.  Consequently, it is crucial to reduce any short-wavelength oscillations in the subdivision surface. As the included examples demonstrate, meshes with extraordinary vertices usually lead to lower convergence rates than meshes with tensor-product connectivity.  For obtaining the improved isogeometric analysis adapted refinement rules  we postulate a constrained optimisation problem with a cost function measuring the errors in approximating cup- and saddle-like quadratic shapes. Three of the weights in the Catmull-Clark subdivision scheme around an extraordinary vertex are chosen as degrees of freedom for optimisation. As constraints the $C^1$ continuity of the surface is strictly enforced and bounded curvatures are enforced as long as non-negative real weights are feasible. The eigenstructure of the subdivision matrix is extensively used in formulating the optimisation problem as usual in previous related work~\cite[Chapter~4,5]{Peters:2008aa} and \cite[Chapter~15]{Sabin:2010aa}. We compute the eigenvalues and eigenvectors numerically after applying a discrete Fourier transform that exploits the local circular symmetry around the extraordinary vertex. The local parameterisation of the subdivision surface required for evaluating the finite element integrals and the cost function is obtained with the algorithm proposed by Stam~\cite{Stam:1998aa}. Two sets of optimised weights for cup- and saddle-like shapes are obtained. The weights for finite element analysis are chosen depending on the dominant shape of the solution field around an extraordinary vertex.

For completeness, we note that subdivision is not the only approach for creating smooth surfaces on arbitrary connectivity control meshes. Over the years numerous $C^k$ and $G^k$ smooth constructions with $k \ge 1$ have been proposed, too many to name here. The search for sufficiently flexible smooth surface representations, especially with  $C^{k \ge 2}$ and $G^{k \ge 2}$, is still open.  
 It is worth mentioning that none of the existing constructions is widely used in commercial CAD systems. This may well be because their implementation is too complicated. The application of basis functions resulting from  smooth  constructions for  isogeometric analysis is currently a very active area of research. For instance, the utility of $G^k$ constructions with NURBS has recently been explored in~\cite{nguyen2016c, collin2016analysis, kapl2016isogeometric}. Alternatively, $C^k$ constructions relying on manifold-based surface constructions~\cite{majeedCirak:2016, ying2004simple,grimm1995modeling} and constructions relying on singular parameterisations have also been investigated~\cite{toshniwal2017smooth, nguyen2016refinable, reif1998turbs}. Some of these schemes are able to provide optimal convergence rates.

The outline of this paper is as follows. In Section~2 the Catmull-Clark subdivision is introduced, with a review of the relevant theory on eigenanalysis of the subdivision matrix. Specifically, the necessary conditions for $C^1$ smoothness and boundedness of the curvature are motivated, and the local parameterisation of subdivision surfaces using the characteristic map is introduced. These are all classical results and concepts which are mostly unknown in isogeometric analysis. In Section~3 the proposed constrained optimisation problem and its numerical solution are discussed. Two sets of subdivision weights are derived that minimise the thin-plate energy norm errors in approximating locally cup- and saddle-like shapes. Subsequently, it is shown how a finite element solution can be locally decomposed into cup- and saddle-like components. Depending on this decomposition and the following choice of optimal weights, a second more accurate finite element analysis can be performed. In Section~4 the proposed approach is applied to transversally loaded thin-plate problems using meshes with extraordinary vertices and the convergence of the errors in $L_2$ and energy norms is reported.

%
\section{Catmull-Clark subdivision surfaces \label{sec: subdivision}}
%
\subsection{Refinement weights and the subdivision matrix}
\begin{figure}[ht!]
	\centering 
  	\subfloat[][Face vertex]
	{
		\includegraphics[scale=0.75]{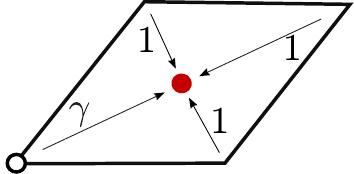}
	}
	\hfil
	\subfloat[][Edge vertex]
	{
		\includegraphics[scale=0.75]{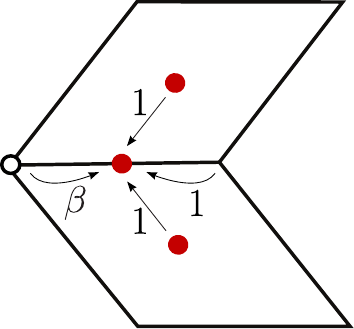}
	}\\	
		\subfloat[][Extraordinary vertex]
	{
		\includegraphics[scale=0.75]{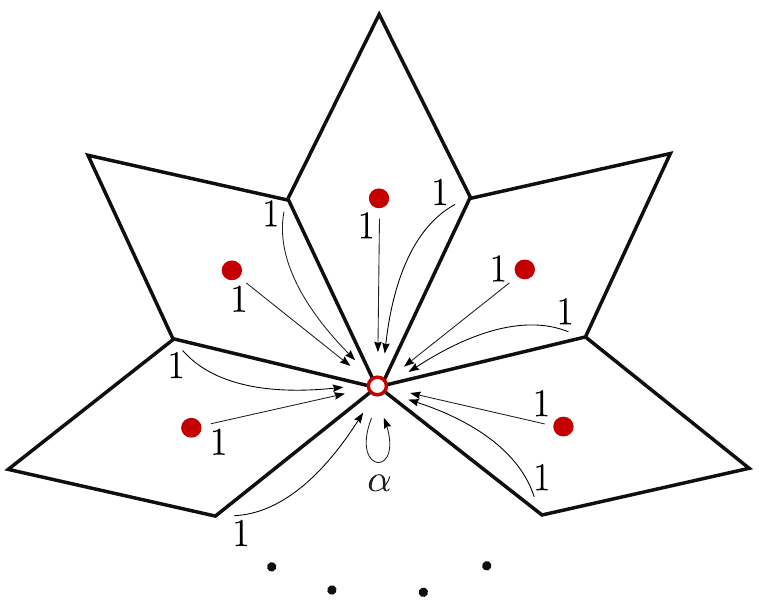}
	}
    	\caption{Subdivision weights for the Catmull-Clark scheme with the empty circle denoting the extraordinary vertex. The weights in each of the three diagrams have to be normalised so that they add up to one. For Catmull-Clark scheme the three weights take the values $\alpha=v(v-2)$, $\beta=1$ and $\gamma=1$, where $v$ is the valence.}
  \label{fig:stencillsCC}
\end{figure}

%
Catmull-Clark subdivision is a generalisation of cubic tensor-product B-splines to unstructured meshes~\cite{Catmull:1978aa}. On non-tensor-product meshes the number of faces connected to a vertex, i.e. valence~$v$, can be different from four. The vertices with $v \neq 4$ are referred to as \textit{extraordinary} or \textit{star vertices}. During subdivision refinement each face of the control mesh is split into four faces and the coordinates of the old and new control vertices are computed with the subdivision weights given in Figure~\ref{fig:stencillsCC}. The weights in each of the three diagrams have to be normalised so that they add up to one. The unnormalised weights assigned to the extraordinary vertex (empty circle) are denoted by $\alpha$, $\beta$ and $\gamma$ respectively. For $v=4$ and bivariate cubic B-splines the three weights take the values $\alpha=8$, $\beta=1$ and $\gamma=1$. The new vertices introduced by the subdivision process are all regular (with  $v = 4$) and the total number of irregular vertices in the mesh remains constant. That is, the irregular vertices are more and more surrounded by regular vertices. 

In order to study the smoothness behaviour of subdivision surfaces near an extraordinary vertex, it is sufficient to consider only the vertices in its immediate vicinity. A $1$-neighbourhood of a vertex is formed by the union of faces that contain the vertex. The $n$-neighbourhood is defined recursively as the union of all $1$-neighbourhoods of the $(n-1)$-neighbourhood vertices.  It is assumed that the considered $n$-neighbourhood has only one single extraordinary vertex located at its centre. The $n$-neighbourhood  control vertices $\vec p^{\ell}$ at the refinement level $\ell$ are mapped to control vertices $\vec p^{\ell+1}$ with the subdivision matrix $\vec S$,  
\begin{equation}
	\vec p^{\ell+1} = \vec S \vec p^\ell \, .
\end{equation}
The square  subdivision matrix $\vec S$ can be readily derived from the weights indicated in Figure~\ref{fig:stencillsCC}. The control point coordinates at level $\ell$ are arranged in this form
\begin{equation}
	\vec p^{\ell} = 
	\begin{bmatrix}
		 p_{1x}^\ell &   p_{1y}^\ell  &  p_{1z}^\ell   \\
		  p_{2x}^\ell &   p_{2y}^\ell &   p_{2z}^\ell    \\
		 \vdots & \vdots & \vdots 
	\end{bmatrix}
\end{equation} 
with each row containing the coordinates of one control point  $\vec p^\ell_j \in \mathbb R^3$ with the index $j$.

%
\subsection{Eigendecomposition of the subdivision matrix \label{sec:eigenDecompos}}
\begin{figure}[ht!]
	\centering 
	\includegraphics[scale=0.5]{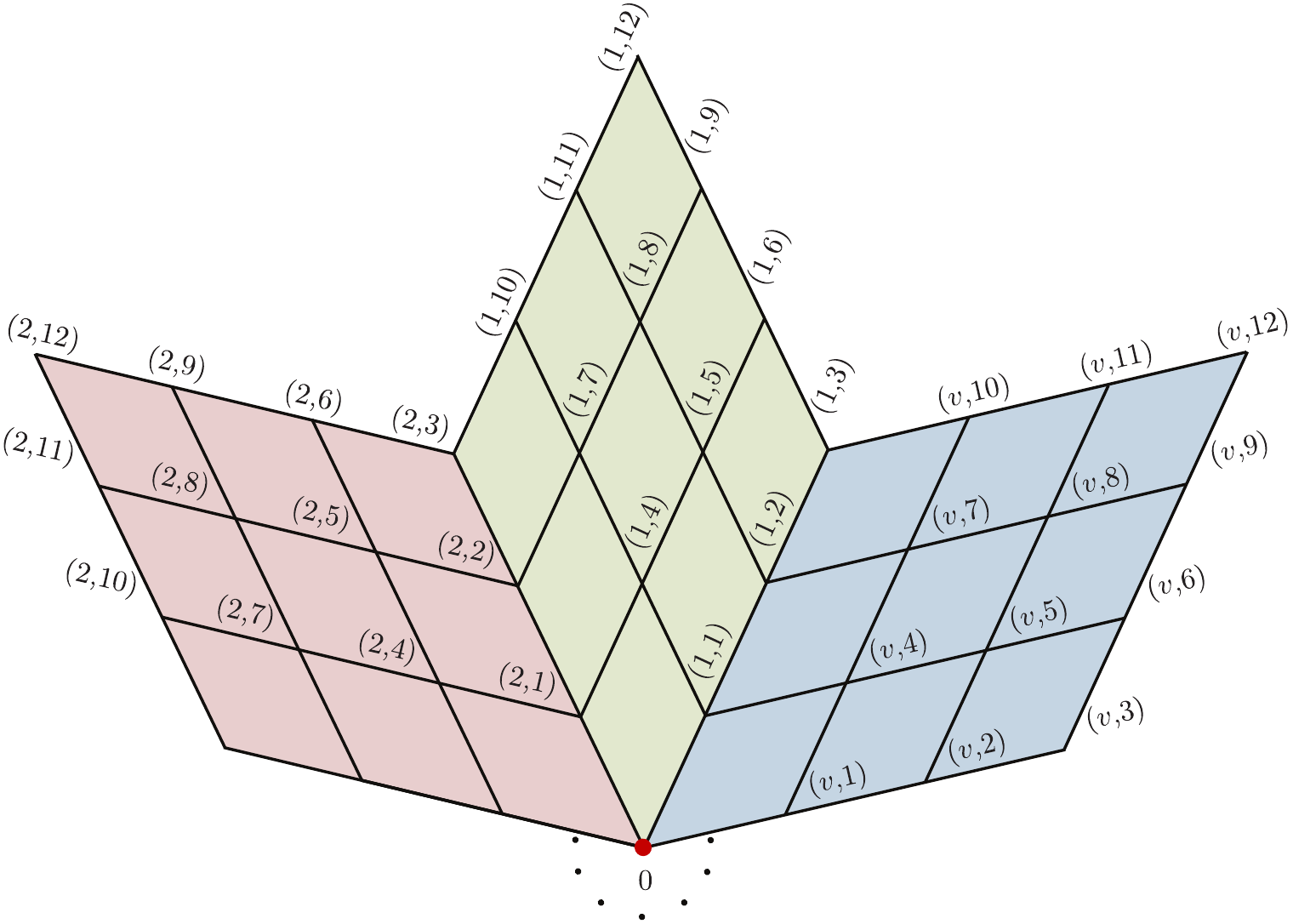}
	\caption{Three-rings of faces around an extraordinary vertex with valence $v$ and the numbering of the vertices. In the index pair $(s,a)$ the first  is the segment number and the second is the vertex number.   \label{fig:oneRingNumbering}}
\end{figure}

%
For the tuning approach to be introduced in Section~\ref{sec: tuning}, it is necessary to consider the $3$-neighbourhood around an extraordinary vertex, see Figure~\ref{fig:oneRingNumbering}.
The  $3$-neighbourhood consists of $v$ segments with each segment containing  $12v$ vertices, excluding the extraordinary vertex with index~$0$. Hence, there are $(12v+1)$ vertices so that the subdivision matrix has the dimensions $(12v+1) \times (12v+1)$.  In establishing the subdivision matrix it is assumed that the index pair  $(s,a)$ is converted to a scalar index as $a+12|s-v|$. The eigenvalues and eigenvectors of $\vec S$ are closely related to the smoothness and other properties of the subdivision surface. The eigendecomposition of the asymmetric subdivision matrix $~\vec S$ reads
\begin{equation}
	\vec S = \sum_j  \lambda_j \vec r_j \otimes \vec l_j  
\end{equation}
with 
\begin{equation} \label{eq:evOrtho}
\begin{split}
	&(\vec S - \lambda_j \vec I) \vec r_j  = \vec 0, \quad  \vec l_j^\trans  (\vec S - \lambda_j \vec I)  = \vec 0^\trans \quad 
	\text { and } \quad 	 \\
	& \langle \vec l_j, \vec r_k \rangle = 
	\begin{cases}
		1 & \text{ if } j = k \\
		0 & \text { if } j \neq k
	\end{cases} \, ,
\end{split}
\end{equation}
where $\lambda_j$ are the eigenvalues and $\vec r_j$, $\vec l_j$ are the right and left eigenvectors respectively. Throughout the paper it is assumed that the eigenvalues are sorted in descending order with largest being $\lambda_0$. The subdivision matrix $\vec S$ has a cyclical structure due to the cyclic symmetry of the weights given in Figure~\ref{fig:stencillsCC}. We assume that the vertices in the $3$-neighbourhood are enumerated according to Figure~\ref{fig:oneRingNumbering}.

Because of the cyclical structure, the eigendecomposition of $\vec S$ can be best computed with a discrete Fourier transform (DFT). As pioneered in~\cite{Doo:1978aa}, DFT is crucial in identifying the different geometric shapes described by the different eigenvectors. For transforming $\vec S$ the following extended DFT matrix is considered:
%
\begin{equation} \label{eq:forwardDFT}
	\vec F  = 
		\frac{1}{\sqrt{v}}
	\begin{bmatrix}
	        1& \vec 0^\trans & \vec 0^\trans & \vec 0^\trans  & \cdots & \vec 0^\trans \\ 
		\vec 0 & \vec I  &  \vec I & \vec I & \cdots & \vec I \\
		\vec 0 & \vec I & {\omega} \vec I & {\omega}^2  \vec I & \cdots &{\omega}^{-1} \vec I \\ 
		\vec 0 & \vec I & {\omega}^2 \vec I & {\omega}^4 \vec I & \cdots & {\omega}^{-2} \vec I \\ 
		\vec 0 & \vdots & \vdots  & \vdots & \ddots & \vdots \\
		\vec 0 & \vec I & {\omega}^{-1} \vec I & {\omega}^{-2} \vec I & \cdots & {\omega} \vec I
	\end{bmatrix}
\end{equation}
with  the complex number ${\omega} = \exp(i 2 \pi/v) $ where $i=\sqrt{-1}$, the identity matrix~$\vec I$ of size $(12 \times 12)$ and  the zero vector~$\vec 0$ of size $12$. The first row and column of~$\vec F$ have been introduced for the extraordinary vertex. In obtaining~\eqref{eq:forwardDFT} the standard relations $\omega^{v+k} = \omega^k $ and $\omega^{v-k} = \omega^{-k} = \overline{\omega}^{k}$ with the complex conjugate $\overline {\omega} = \exp(-i 2 \pi/v)$ are used. The inverse transform $\vec F^{-1}$ is obtained by replacing $\omega$ with its complex conjugate $\overline {\omega}$. The subdivision matrix is Fourier transformed according to
\begin{equation}
	\hat {\vec S} = \vec F \vec S \vec F^{-1} \, ,
\end{equation}
leading to a block diagonal  matrix
\begin{equation}
	\hat {\vec S} = 
	\begin{bmatrix}
		\hat {\vec S}^{(0,0)} &  & &  & \\
		& \hat {\vec S}^{(1,1)}   & &  & \\ 
		& & \hat {\vec S}^{ (2,2)}   & & \\  
		& & & \ddots & \\
		& & &  & \hat {\vec S}^{ (v-1, v-1)}    \\  		
	\end{bmatrix} \, ,
\end{equation}
where the blocks $\hat {\vec S}^{(m,m)}$ are of size $13\times13$ for $m=0$ and of size $12\times12$ for $m \neq 0$.  Due to its block-diagonal structure the eigendecomposition of the transformed matrix  
\begin{equation}
	\hat{\vec S} = \sum_j \lambda_j \hat {\vec r}_j \otimes \hat {\vec l}_j
\end{equation}
can be more readily determined. Namely, it is sufficient to consider the eigenvalue problems for each of the $v$ blocks $\hat{\vec S}^{(m,m)}$ separately, i.e., 
\begin{equation}
	\left ( \hat {\vec S}^{(m,m)}  - \lambda^{(m,m)}_n \vec I \right )  \hat {\vec r}^{(m,m)}_n = \vec 0 \, ,
\end{equation}
where the eigenvalues within each block are also sorted in descending order, i.e, $\lambda^{(m,m)}_0$ is the largest eigenvalue in the block $\hat {\vec S}^{(m,m)}$. The eigenvalues $\lambda_j$ and eigenvectors $\hat {\vec r}_j$ and $\hat {\vec l}_j$ are the union of all the eigenvalues and the block-wise eigenvectors. In obtaining the eigenvectors  $\hat {\vec r}_j$ and $\hat {\vec l}_j$, each of size $12 v+1$,  the corresponding block-wise vectors $\hat {\vec r}^{(m,m)}_n$  and $\hat {\vec l}^{(m,m)}_n$   are  suitably padded with zeros.  The vectors $\hat {\vec r}^{(m,m)}_n$  and $\hat {\vec l}^{(m,m)}_n$ are of size $13$ for $m=0$ and of size $12$ for $m \neq 0$.
Moreover, the subdivision matrix ${\vec S}$ and its Fourier transform $\hat{\vec S}$ have the same eigenvalues $\lambda_j$ and their eigenvectors are related by 
\begin{equation} \label{eq:inverseTransform}
 	\vec r_j = \vec F^{-1} \hat {\vec r}_j \quad \text{ and } \quad \vec l_j = \vec F^{-1} \hat { \vec l}_j \, . 
\end{equation}

Each block $\hat {\vec S}^{(m,m)}$ corresponds to a specific rotational frequency $ \omega_f = 2 \pi m/ v$. As pointed out, the eigenvectors  $\hat {\vec r}_j$ and $\hat {\vec l}_j$  can have non-zero entries only in the components corresponding to a specific $\hat {\vec r}^{(m,m)}_n$  and $\hat {\vec l}^{(m,m)}_n$. Hence, the transformation of  $\hat {\vec r}_j$ and $\hat {\vec l}_j$ according to~\eqref{eq:inverseTransform} yields always a  column of $\vec F^{-1}$ each of which corresponds to a specific rotational frequency\footnote{The first columns and rows of $\vec F$ and $\vec F^{-1}$ are assigned to the extraordinary vertex and do not represent harmonics.}.  To this end, recall the Euler identity 
\begin{equation}
	\omega^{ms} =  e^{i 2 \pi  ms/v} = \cos (2 \pi ms /v) + i \sin (2 \pi ms / v) \, . 
\end{equation}
Hence, for a fixed angular frequency $\omega_f = 2 \pi m/ v$  the vectors $\vec r_j$ and $\vec l_j$  will assign each control vertex $(s,a)$ with a fixed index $a$, c.f. Figure~\ref{fig:oneRingNumbering}, a value that oscillates with the angular frequency $2 \pi m /v$ while circumnavigating the extraordinary vertex by incrementing $s \in \{1, \cdots, v \}$.


Furthermore, for geometric interpretation of  the eigendecomposition it is helpful to realise that most of the eigenvalues $\lambda_j$ have the multiplicity of two. That is, the eigenvalues  $ \lambda_n^{(m,m)}$ and $ \lambda_n^{(v- m, v-m)}$ are identical in the blocks $m \ge 1$. The corresponding eigenvectors $\vec r_j$ and $\vec l_j$  have the same eigenfrequency because the columns $m$ and $v-m$ of the DFT matrix~$\vec F^{-1}$  are the complex conjugates of each other. 

%
\subsection{Limit analysis and smoothness \label{sec:LimitSmooth}}
%
The eigenvalues $\lambda_j$ and eigenvectors $\hat{\vec r}_j$ and $\hat {\vec l}_j $ of the Fourier transformed subdivision matrix $\hat{\vec S}$ have to satisfy certain conditions for a subdivision scheme leading to a smooth well-defined surface, see~\cite{Doo:1978aa,Peters:2008aa}. To understand this, consider the projection of control mesh vertex coordinates  at subdivision level $\ell=0$  into the eigenspace of the subdivision matrix, using the orthogonality of left and right eigenvectors~\eqref{eq:evOrtho}, 
\begin{equation} \label{eq:cpProjectEigen}
	\begin{split}
	\vec p ^{\ell=0} = & ~\vec r_0 \langle \vec l_0 , \vec p^0 \rangle +  \vec r_1 \langle \vec l_1 , \vec p^0 \rangle +  \vec r_2 \langle \vec l_2 , \vec p^0 \rangle + \cdots   \\
	& +  \vec r_{12v} \langle \vec l_{12v} , \vec p^0 \rangle  \, ,
	\end{split}
\end{equation}
where each of the scalar products $ \langle \phantom{x},\phantom{x} \rangle$ yield a row vector with $3$ components. Subdividing the  $3$-neighbourhood in the eigenspace, while considering the eigendecomposition~\eqref{eq:evOrtho},   gives
\begin{equation}
	\begin{split}
		\vec S \vec p ^{0}  = & ~\lambda_0 \vec r_0 \langle \vec l_0 , \vec p^0 \rangle +  \lambda_1  \vec r_1 \langle \vec l_1 , \vec p^0 \rangle +  \lambda_2  \vec r_2 \langle \vec l_2 , \vec p^0 \rangle + \cdots \\
		& + \lambda_{12v}  \vec r_{12v} \langle \vec l_{12v} , \vec p^0 \rangle  \, . 
	\end{split}
\end{equation}
Hence, the repeated subdivision of the  $3$-neighbourhood can be simply achieved with
\begin{equation} \label{eq:limit3Ring}
	\begin{split}
		\vec S^\ell \vec p ^{0}  = & ~\lambda_0^\ell \vec r_0 \langle \vec l_0 , \vec p^0 \rangle +  \lambda_1^\ell  \vec r_1 \langle \vec l_1 , \vec p^0 \rangle +  \lambda_2^\ell  \vec r_2 \langle \vec l_2 , \vec p^0 \rangle + \cdots   \\
 & + \lambda_{12v}^\ell  \vec r_{12v} \langle \vec l_{12v} , \vec p^0 \rangle  \, . 
	\end{split}
\end{equation}
From this equation it is evident that the properties of a subdivision surface are widely governed by the eigenstructure of the subdivision matrix. The subdivision matrix $\vec S$ is a stochastic matrix, i.e. only positive entries and each row adds up to 1, so that its largest eigenvalue is $\lambda_0 = 1$ and the components of the corresponding eigenvector $\vec r_0$ are all equal to~$1$.  In the limit $\ell \rightarrow \infty$ all control vertices converge to $\langle \vec l_0 , \vec p^0 \rangle$. The first term in~\eqref{eq:limit3Ring} can be eliminated by translating the initial control vertex coordinates by $- \vec r_0 \langle \vec l_0, \vec p^0 \rangle$. Without loss of generality, in the following we assume that the coordinate system for $3$-neighbourhood has been chosen so that the first term in~\eqref{eq:limit3Ring} is zero, that is, 
\begin{equation} \label{eq:limit3RingT}
\begin{split}
		\vec S^\ell \vec p ^{0}  =& ~ \lambda_1^\ell  \vec r_1 \langle \vec l_1 , \vec p^0 \rangle +  \lambda_2^\ell  \vec r_2 \langle \vec l_2 , \vec p^0 \rangle  + \cdots 
		+ \lambda_{12v}^\ell  \vec r_{12v} \langle \vec l_{12v} , \vec p^0 \rangle  \, . 
\end{split}
\end{equation}
For a (symmetric)  $C^1$-continuous subdivision surface the subdominant eigenvalues $\lambda_1$ and $\lambda_2$ have to satisfy the following relationship:
\begin{equation} \label{eq:c1constraints}
	\lambda_1 = \lambda_2 > \lambda_3  \, .
\end{equation}
In addition, the corresponding eigenvectors $\vec r_1$, $\vec l_1$, $\vec r_2$ and $\vec l_2$ have to come from the eigendecomposition of the blocks $\vec {\hat S}^{(1,1)}$ and $\vec {\hat S}^{(v-1,v-1)}$~\cite{Doo:1978aa, Reif:1995aa}. As discussed in Section~\ref{sec:eigenDecompos},  owing to the symmetry properties of the Fourier transformation, $\vec {\hat S}^{(1,1)}$ and $\vec {\hat S}^{(v-1,v-1)}$ have the same eigenvalues, and the eigenvectors    $\vec r_2$ and $\vec l_2$   are the complex conjugates of $\vec r_1$ and $\vec l_1$. All the four eigenvectors $\vec r_1$, $\vec l_1$, $\vec r_2$ and $\vec l_2$ are usually complex and have  the angular frequency $\omega_f=2 \pi/v$.  A set of real eigenvectors each of size $12v +1$ representing vertex values can be obtained as the linear combination of the complex ones, e.g., with   $ \tfrac{1}{2} (\vec r_1 + \vec r_2) $ and $  \tfrac{1}{2i} (\vec r_1 - \vec r_2)$ where $i=\sqrt{-1}$. To avoid a proliferation of symbols we will use the same symbols for the so-computed real and complex eigenvectors. 

A necessary condition for the $C^2$-continuity of a subdivision surface (with no artificial flat spots) is that the subsubdominant eigenvalues satisfy 
\begin{equation}
	\label{eq:c2constraints}
	\lambda_3 = \lambda_1^2,  \quad  
	\lambda_4 = \lambda_1^2 \, , \quad \lambda_3 = \lambda_4 = \lambda_5 > \lambda_6
\end{equation}
and the corresponding eigenvectors come from the eigendecomposition of the blocks $\vec {\hat S}^{(0,0)}$, $\vec {\hat S}^{(2,2)}$ and $\vec {\hat S}^{(v-2,v-2)}$~\cite{Peters:2008aa}. Remember that $\lambda_4 = \lambda_5$ is naturally satisfied due to the duplicity of eigenvalues from blocks $\vec {\hat S}^{(m,m)}$ and $\vec {\hat S}^{(v-m,v-m)}$ when $m \ge 1$, as mentioned in Section~\ref{sec:eigenDecompos}. 
%
\subsection{Characteristic map}
\label{sec:characteristicMap}
\begin{figure}[ht!]
	\centering 
	\includegraphics[width=0.18\textwidth]{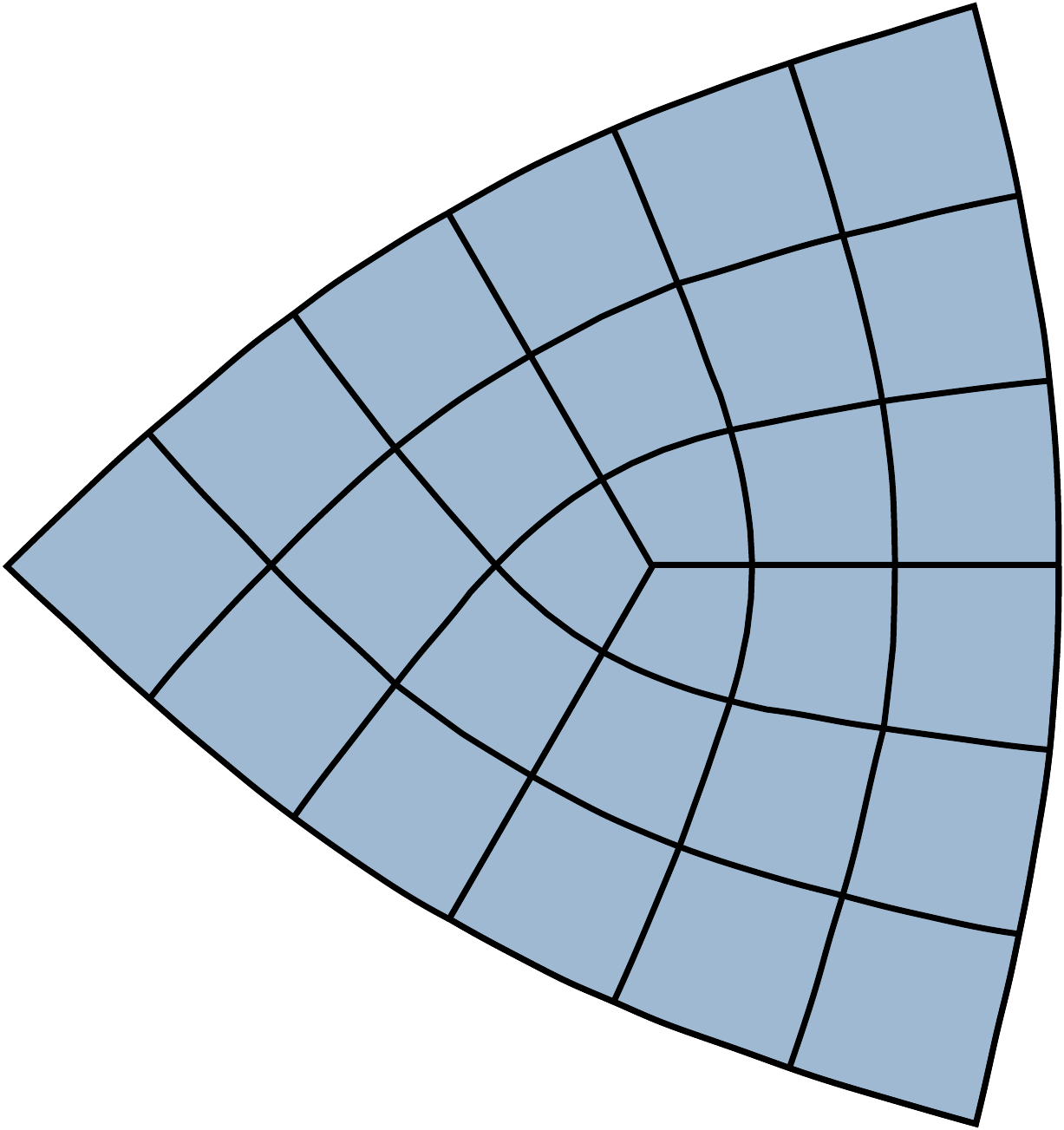}
	\hfil
	\includegraphics[width=0.18\textwidth]{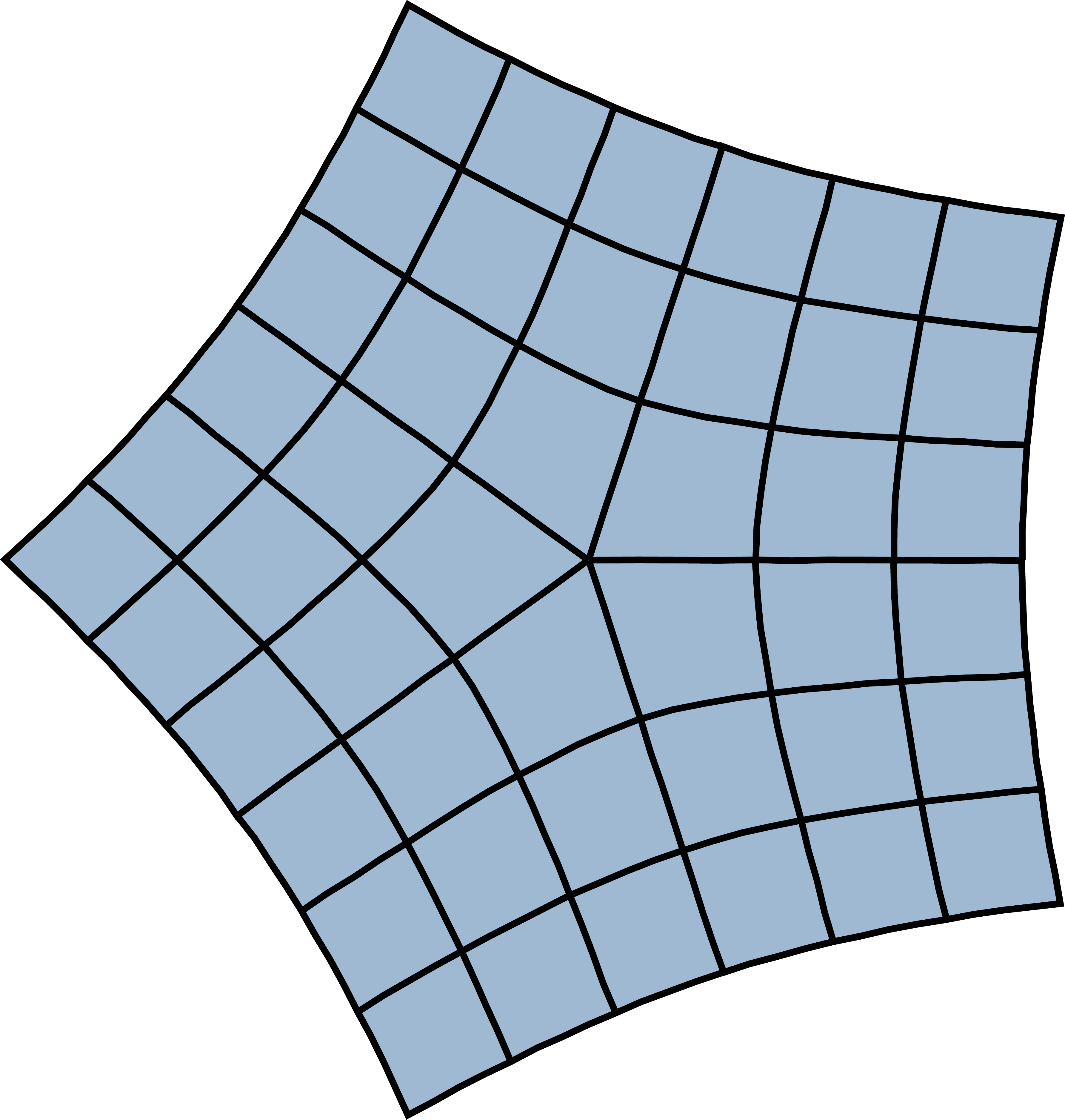}
	\caption{Characteristic control mesh of Catmull-Clark scheme for valence $v=3$ (left) and $v=5$ (right). 
	\label{fig:characteristicMesh}}
\end{figure}

%
As first proposed in Reif~\cite{Reif:1995aa} the characteristic map provides a means for parameterisation of the surface generated by a subdivision scheme. Parameterisation of the subdivision surface, at least a local one, is essential in order to associate the so far discrete representation based on control vertices with  a continuous differentiable representation. This is, for instance, required for finite element analysis using subdivision surfaces. The characteristic map is defined using the two real right eigenvectors $\vec r_1$  and $\vec r_2$ corresponding to the subdominant eigenvalue $\lambda_1 = \lambda_2$. The characteristic control mesh shown in Figure~\ref{fig:characteristicMesh} representing the 3-neighbourhood around an extraordinary vertex has the coordinates 
\begin{equation}
	\vec p^{0}_c = 
	\begin{bmatrix}
		 \vec r_{1} &    \vec r_{2}    & \vec 0    
	\end{bmatrix} \, ,
\end{equation} 
where  the third \emph{out-of-plane} coordinate is chosen as $\vec 0$.  As discussed in Section~\ref{sec:eigenDecompos}, recall that the two eigenvectors $ \vec r_{1}$ and $ \vec r_{2}$ have the angular frequency $\omega_f = 2 \pi/v$ and have been chosen so that they are orthogonal in the plane spanned by the corresponding two complex eigenvectors. This can be done without loss of generality because the subdivision construction is invariant under affine transformations. Hence, the coordinates of  control vertices $(s,a)$ oscillate with $\cos (2\pi s /v) $ in the horizontal direction and with $\sin (2 \pi s/v) $ in the vertical direction leading to the shown characteristic control mesh.  
As suggested in Reif~\cite{Reif:1995aa}, the planar surface described by the characteristic control mesh can be used for the parameterisation of subdivision surfaces. To this end, first consider the subdivision refinement of the characteristic mesh. According to~\eqref{eq:limit3RingT} and the orthogonality of left and right eigenvectors~\eqref{eq:evOrtho}, the subdivision refinement of the characteristic mesh simply yields a scaled version of the same mesh: 
\begin{equation}
	\vec p^{\ell}_c  = \vec S^\ell \vec p^{0}_c = 
	\begin{bmatrix}
		\lambda_1^\ell \vec r_{1} &  \lambda_2^\ell  \vec r_{2}    & \vec 0     
	\end{bmatrix} \, .
\end{equation}  
Hence, the refined control mesh is simply obtained by scaling the control mesh by $\lambda_1=\lambda_2$. Repeated subdivision yields repeated scaling of the control mesh. This combined with the fact that the Catmull-Clark scheme leads to bivariate cubic B-splines in patches with only ordinary vertices is used for parameterising the subdivision surface. During subdivision refinement each patch is split into four patches. In particular, in the patches adjacent to the extraordinary vertex three of the created patches have only regular vertices and can be parameterised with bivariate cubic B-splines. 

\begin{figure*}[ht!]
	\centering 
	\includegraphics[width=0.58\textwidth]{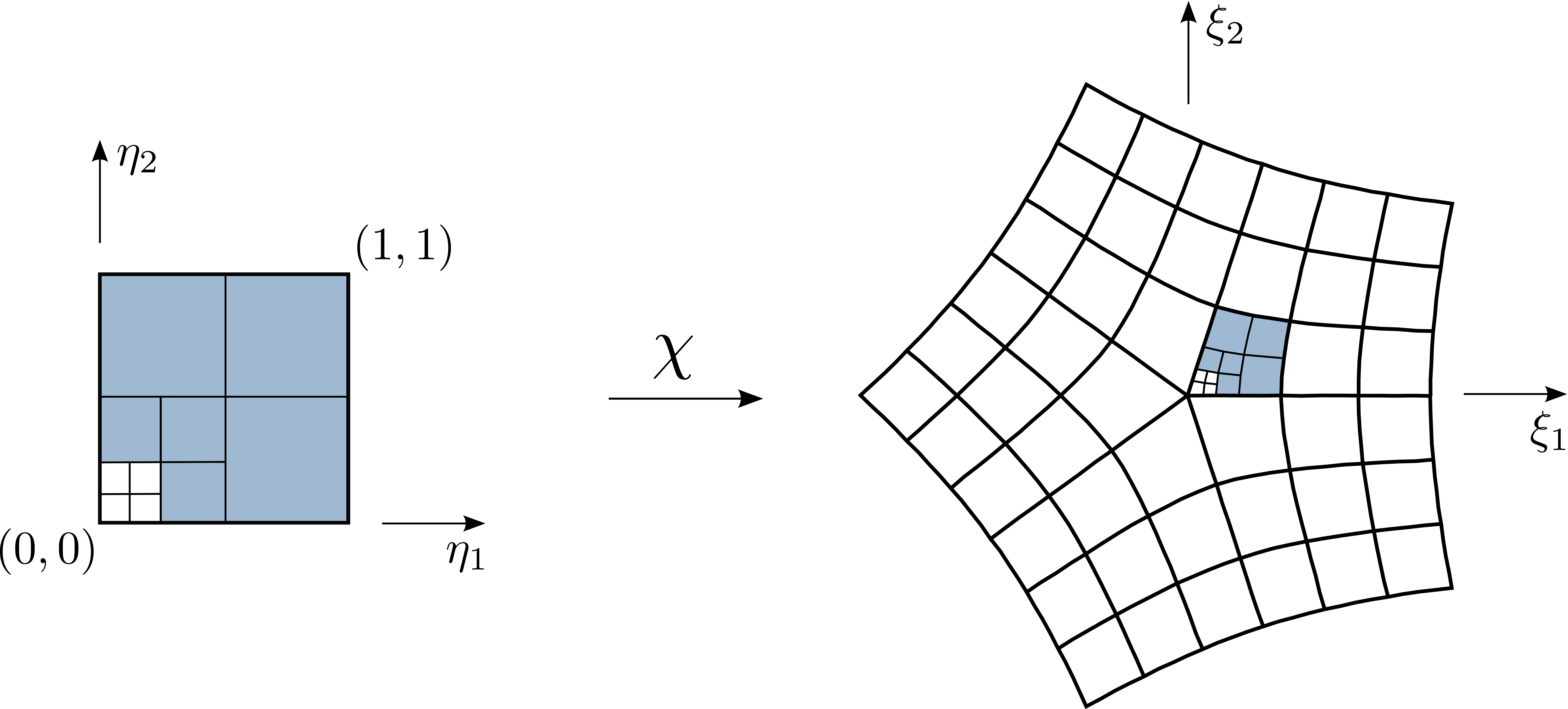}
	\caption{Characteristic map from a unit square to the characteristic domain. 
	\label{fig:characteristicMap}}
\end{figure*}

With repeated refinement more and more of the subdivision surface can be parameterised with cubic B-splines. A practical algorithm for efficient implementation of this parameterisation has been introduced in Stam~\cite{Stam:1998aa}. Without going into details we define the bijective characteristic map 
%
%
\begin{equation}
	 \chi :  (\vec \eta , s)  \in (\Omega, s)  \mapsto  \vec \xi \in  \Omega_\chi  \qquad 
\end{equation}
with $\vec \eta = ( \eta_1, \eta_2  )$ and  $ \vec \xi=(\xi_1, \xi_2)$,  
which maps  a set of square domains $(\Omega, s)$ with $s \in \mathbb N^+$ representing the faces in the control mesh into the characteristic domain $\Omega_\chi$.  The smooth parameterisation provided by the characteristic map $\chi$ is illustrated in Figure~\ref{fig:characteristicMap}.
With the subdivision basis functions $\vec N (\vec \eta, s) $, consisting of cubic B-splines and obtained according to~\cite{Stam:1998aa}, the characteristic map can be written as 
\begin{equation} \label{eq:charactInterp}
	\vec \xi = \chi (\vec \eta, s)  = {\vec N}^\trans (\vec \eta, s )  
		\begin{bmatrix}
		\vec r_{1} &  \ \vec r_{2}        
	\end{bmatrix} \, .
\end{equation}
For brevity, in  the following we omit  the face index in the basis function ${\vec N} (\vec \eta, s )$  and write ${\vec N} (\vec \eta )$.

%
\section{Optimisation of subdivision weights \label{sec: tuning}}
%
We aim to modify the subdivision weights $\alpha, \beta$ and $\gamma$ of the Catmull-Clark scheme to improve its approximation properties when used in finite element analysis. As known in CAD, not all parameters yield visually appealing surfaces even when they give $C^1$-continuous surfaces, for a quantitative analysis see~\cite{peters2004shape}. Small surface oscillations, i.e. ripples, appear when the represented surface is not planar. 

\subsection{Preliminaries}
\label{sec:preliminaries}
%

First, we consider the representation of a polynomial scalar field $u(\xi_1, \xi_2 )$ over the characteristic domain $\Omega_\chi$. It is assumed that the scalar field is given in the form
\begin{equation} \label{eq:wSeries}
	\begin{split}
	u(\xi_1, \xi_2 )  = & ~c_0  + c_1 \xi_1 + c_2 \xi_2 + c_3 (\xi_1^2+ \xi_2^2)  + c_4 (\xi_1^2-\xi_2^2) \\
	& + c_5 (2 \xi_1 \xi_2) + \ldots  	 \\ 
	= & ~c_0 u_{0} +  c_1 u_{1} (\xi_1) +  c_2 u_{2}  (\xi_2)+   c_3 u_{3} (\xi_1, \xi_2) \\ 
	& + c_4 u_{4} (\xi_1, \xi_2)  +  c_5 u_{5} (\xi_1, \xi_2) +  \ldots  \, ,
	\end{split}   
\end{equation}
where $c_j \in  \mathbb R$  and the  functions $u_j$ on the second line are introduced for notational convenience. The chosen functions $u_0$, $u_1(\xi_1)$, $u_2(\xi_2)$, $u_3(\xi_1, \xi_2)$, $u_4(\xi_1, \xi_2)$ and $u_5(\xi_1, \xi_2)$ can represent all quadratics and their choice will be discussed further below. The approximation of $u_j$  over $\Omega_\chi$ can be studied by comparing with it the limit surface resulted from the control points
\begin{equation} \label{eq:controlCupHyp0}
	\vec p^{0}_{u_j} = 
	\begin{bmatrix}
		 \vec r_{1} &    \vec r_{2}    &  u_j(\vec r_{1}, \vec r_{2} ) 
	\end{bmatrix} \, , 
\end{equation}   
where the third coordinate is a vector formed by the scalar function $u_j(\xi_1, \xi_2 )$ evaluated at the vertex locations $[ \vec r_1 \; \vec r_2]$, row by row. The linear functions~$u_1(\xi_1)$ and $u_2(\xi_2)$ can be exactly represented so that we are mainly concerned about the quadratic terms $u_3(\xi_1, \xi_2)$, $u_4(\xi_1, \xi_2)$ and $u_5(\xi_1, \xi_2)$.

\begin{figure}[ht!]
  \centering
   \subfloat[][Cup-like geometry~$u_3^h (\xi_1, \xi_2)$] 
  {
     \includegraphics[width=0.26\textwidth]{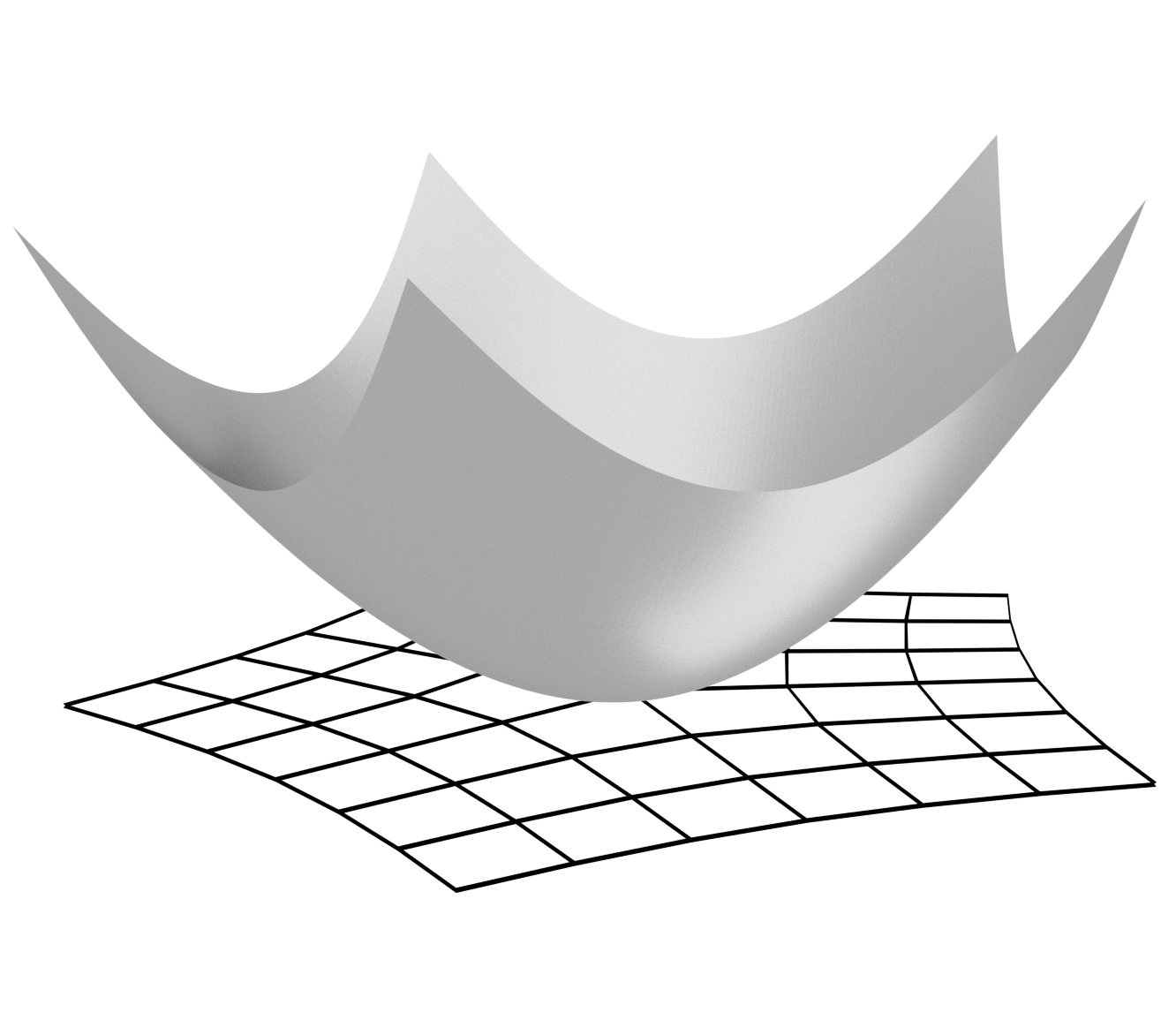}
  } 
  \\
   \subfloat[][Saddle-like  geometry~$u_4^h (\xi_1, \xi_2)$] 
  {
     \includegraphics[width=0.26\textwidth]{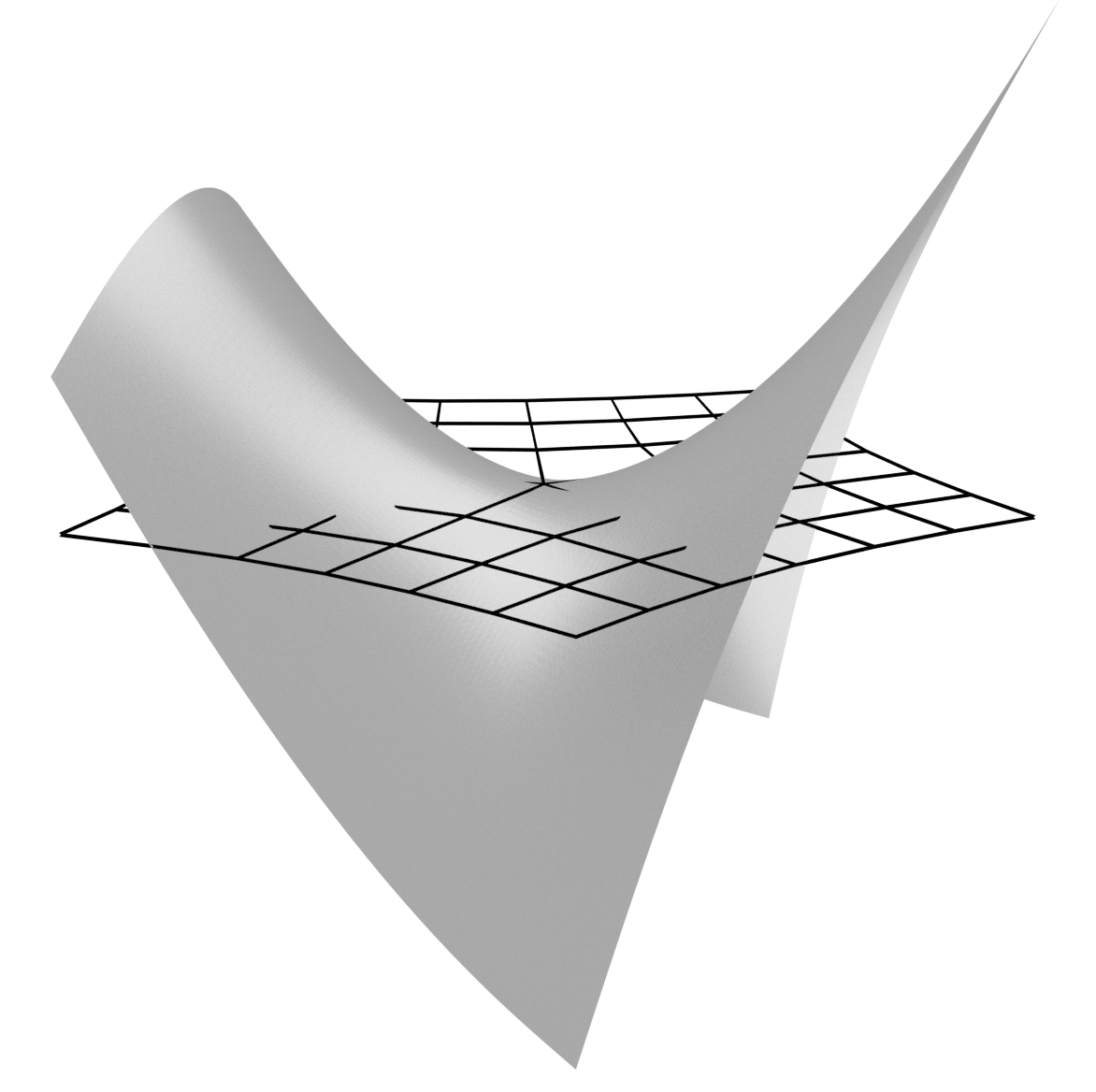}
  }
  \caption{Quadratic shapes over a characteristic control mesh with valence $v=5$. Note that~$u_5^h (\xi_1, \xi_2)$ has the same shape like~$u_4^h (\xi_1, \xi_2)$, but only rotated by $\pi/4$ in the $\xi_1 \xi_2$-plane.  }
  \label{fig:cupAndSaddle}
\end{figure}

The specific form of the quadratic functions in~\eqref{eq:wSeries} is motivated by the eigenstructure of the subdivision matrix~$\vec S$, see Section~\ref{sec:eigenDecompos}. Specifically, the control point values $\vec r_3$, $\vec r_4$ and $\vec r_5$ and the corresponding control point values $u_{3}( \vec r_1, \vec r_2) $, $u_{4}( \vec r_1, \vec r_2) $ and $ u_{5}( \vec r_1, \vec r_2) $ have matching angular frequencies over the 3-neighbourhood of the extraordinary vertex\footnote{In order for the phase to match, the indexing of the vertices has to begin along the edge aligned with the $\xi_1$-axis.}. It is straightforward to confirm the orthogonality relations 
\begin{equation} \label{eq:orthogonality}
	\langle u_{j}(\vec r_1, \vec r_2), \vec l_k \rangle = 
	\begin{cases} 
		\neq 0  & \text { for } j=k \\
		= 0 & \text { for } j \neq k 
	\end{cases}
	\quad  \text {with } j,k \in \{3,4,5\} \, .
\end{equation}
According to~\eqref{eq:cpProjectEigen}, the projection of the control vertex coordinates $\vec p^{0}_{u_j}$ into the eigenspace of the subdivision matrix $\vec S$, while neglecting the terms with higher orders than quadratic, yields
\begin{equation} \label{eq:controlCupHyp1}
	\vec p^{0}_{u_j} = 
	\begin{bmatrix}
		 \vec r_{1} &    \vec r_{2}    &  \vec {r}_j \langle \vec l_j,  u_j(\vec r_{1}, \vec r_{2}) \rangle  
	\end{bmatrix} \, . 
\end{equation}   
With the eigendecomposition~\eqref{eq:evOrtho} the subdivision refinement of this control mesh gives
\begin{equation} \label{eq:controlCupHypRefined}
	\vec p^{\ell}_{u_j}  = \vec S^\ell  \vec p^{0}_{u_j} = 
	\begin{bmatrix}
		 \lambda_1^\ell \vec r_{1} &   \lambda_2^\ell  \vec r_{2}    &  \lambda_j^\ell \vec r_{j}   \langle \vec l_j,  u_j(\vec r_{1}, \vec r_{2}) \rangle 
	\end{bmatrix} \, , 
\end{equation}   
That is, the subdivision refinement of the first two components yields the characteristic domain and the third component yields the graph of the surface $u_j^h(\xi_1, \xi_2)$ approximating  $u_j(\xi_1, \xi_2)$, see Figure~\ref{fig:cupAndSaddle}.  The corresponding limit surface $u_j^h(\xi_1, \xi_2)$ has, according to~\eqref{eq:charactInterp}, the following form: 
\begin{equation}
	\begin{bmatrix}
		\vec \xi &  u_j^h  
	\end{bmatrix}
	= \vec N^\trans  (\vec \eta) 
	\begin{bmatrix}
		\vec r_1 & \vec r_2 & \vec r_{j} \langle \vec l_j, u_j(\vec r_{1}, \vec r_{2}) \rangle
	\end{bmatrix} \, .
\end{equation}
To compare the shapes $u_j^h (\xi_1, \xi_2)$ and $u_j (\xi_1, \xi_2)$ quantitatively, we introduce the  thin-plate energy norm  
\begin{equation} \label{eq:thinPlateEnergy}
	\| u \|_e^2 =  \int_\Omega \left ( \frac{\partial^2 u }{\partial \xi_1^2} +   \frac{\partial^2 u}{\partial \xi_2^2 }  \right  )^2 - 2(1-\mu) \left (  \frac{\partial^2 u}{\partial \xi_1^2 }    \frac{\partial^2 u}{\partial \xi_2^2 }  -  \left ( \frac{\partial^2 u}{\partial \xi_1 \partial \xi_2}  \right )^2  \right)  \D\Omega \, , 
\end{equation}
with the Poisson ratio $\mu = 0.3$. 
Moreover, the necessary conditions for $C^2$-continuity given in~\eqref{eq:c2constraints}, repeated here for convenience, 
\begin{equation*}
	\lambda_3 = \lambda_1^2, \quad \lambda_4 = \lambda_1^2, \quad \text{ and } \quad  \lambda_5 = \lambda_1^2  	 \quad \text{ with } \lambda_1 = \lambda_2
\end{equation*}
can now be related to the curvature of the three quadratic limit surfaces resulted from repeated refinement of $\langle  u_j(\vec r_{1}, \vec r_{2}), \vec l_j \rangle \vec r_{j}$ with $j \in \{3,4,5\}$.  In order for the limit surfaces $u_j^h(\xi_1, \xi_2)$ to have finite curvature at the extraordinary vertex,  when the first two control vertex components scale with $\lambda_1 (=\lambda_2)$ the third has to scale with  $\lambda_1^2$.
 
\subsection{Constrained optimisation}
\label{sec:cupAndSaddleTuning}
%
The constrained optimisation problem for determining the subdivision weights $\alpha$, $\beta$ and $\gamma$ that minimise the error in approximating  quadratic surfaces is formulated as
\begin{subequations}
\begin{alignat}{2} \label{eq:optimisation}
  	 &   \underset{\alpha, \beta, \gamma}{\text{minimise}}  & \quad  & \frac{\|  u_j^h (\xi_1, \xi_2; \alpha, \beta, \gamma)   - u_j (\xi_1, \xi_2) \|_e }{ \| u_j (\xi_1, \xi_2) \|_e } \\
 	 &  \text{subject to:} &  & \lambda_1(\beta, \gamma)  = \lambda_2 (\beta, \gamma) \label{eq:optCstr1} \\
	 &			& & \lambda_3 (\alpha, \beta, \gamma) = 	 \lambda_1^2(\beta, \gamma)  \label{eq:optCstr2} \\ 
  	 &  			&  & \lambda_4 ( \beta, \gamma) = 	 \lambda_1^2(\beta, \gamma)  \label{eq:optCstr3} \\ 
	 & 	 	        & & \lambda_5 ( \beta, \gamma) = 	 \lambda_4( \beta, \gamma)  \label{eq:optCstr4} \, , 
\end{alignat}  
\end{subequations} 
with $j \in \{ 3,4,5\}$ and  the constraints representing the necessary $C^2$-continuity conditions \eqref{eq:c1constraints} and \eqref{eq:c2constraints}. As mentioned in Section~\ref{sec:LimitSmooth}, owing to the symmetries of the DFT, the constraints~\eqref{eq:optCstr1} and~\eqref{eq:optCstr4} are automatically satisfied. Hence, the constraints reduce to two independent equations for the three unknowns. To reduce the constrained optimisation problem into an unconstrained one, it is convenient to first solve  the nonlinear system of equations 
\begin{subequations} \label{eq:constraints1}
	\begin{align}
		\lambda_1(\beta, \gamma) &= \lambda \, , \\
		\lambda_4(\beta, \gamma) &= \lambda_1^2(\beta,\gamma)  \, ,  \label{eq:constraints12}\\
			\lambda_3(\alpha,\beta,\gamma) &= \lambda_1^2(\beta,\gamma)  \, .
	\end{align}
\end{subequations} 
That is, to determine the dependence of the weights~$\alpha(\lambda)$, $\beta(\lambda)$ and~$\gamma(\lambda)$ on the variable~$\lambda$. 
To solve \eqref{eq:constraints1}  we use in our implementation the Python library SciPy, to be more specific, the quasi-Newton method with a BFGS update with a suitable cost function.
However,~$\beta(\lambda)$ and~$\gamma(\lambda)$ can become complex for some~$\lambda$ values~\cite{augsdorfer2006tuning}. For Catmull-Clark, it is smaller~$\lambda$ values which result in complex weights. For instance, there is no real solution for ~$\beta$ and~$\gamma$  for $\lambda \leq 0.608$ in case of valence $v=5$. Instead of excluding~$\lambda$ values leading to complex weights we  relax the second constraint~\eqref{eq:constraints12} by considering the modified constraint equations
\begin{subequations} \label{eq:constraints2}
	\begin{align}
		\lambda_1(\beta, \gamma) &= \lambda \, , \\
		\beta &= \gamma  \, ,  \label{eq:constraints22}\\
			\lambda_3(\alpha,\beta,\gamma) &= \lambda_1^2(\beta,\gamma)  \, .
	\end{align}
\end{subequations} 
In implementations where the boundedness of curvature must be satisfied, one can constrain the $\lambda$ value in a valid range or consider more degrees of freedom for optimisation~\cite{augsdorfer2006tuning, it:2016-012}. In our numerical experiments, we found that considering the modified constraint equations leads to smaller energy norm errors in comparison to constraining the range of possible $\lambda$ values. 
After solving~\eqref{eq:constraints1} or ~\eqref{eq:constraints2} and determining  $\alpha(\lambda)$, $\beta(\lambda)$, and $\gamma(\lambda)$ the constrained optimisation problem~\eqref{eq:optimisation} can now be restated as an unconstrained problem 
\begin{equation} \label{eq:optimisationUC}
  	  \underset{\lambda}{\text{minimise}}  \quad  \frac{\|  u_j^h (\xi_1, \xi_2; \alpha(\lambda), \beta(\lambda), \gamma(\lambda))  - u_j (\xi_1, \xi_2) \|_e }{ \| u_j (\xi_1, \xi_2) \|_e }  \, ,
 \end{equation} 
which is a one-dimensional optimisation problem that can be solved by direct search.

\subsection{Optimised weights for valence $v=5$ vertices}
\label{}
%

As an example for obtaining optimised weights, we consider  the valence $v=5$ vertex case. The proposed optimisation follows the same procedure regardless of valence. It is sufficient to consider only the approximation of the quadratic functions $u_3(\xi_1, \xi_2)= \xi_1^2+\xi_2^2$ and $u_4(\xi_1, \xi_2) = \xi_1^2 - \xi_2^2$ with cup-like and saddle-like geometries, respectively.  The function $u_5(\xi_1, \xi_2) = 2 \xi_1 \xi_2$ has the same saddle-like geometry as  $u_4(\xi_1, \xi_2) $, only rotated by~$\pi/4$ in the~$\xi_1\xi_2$-plane.  During optimisation the thin-plate energy norms in~\eqref{eq:optimisationUC} are evaluated in the 2-neighbourhood  of the extraordinary vertex, which is the same as the support size of the basis functions. This explains why we consider $3$-neighbourhood around an extraordinary vertex, see Figure~\ref{fig:oneRingNumbering}, because the evaluation in the second-ring elements needs the third-ring control vertices. 

\begin{figure}[ht!]
  \centering
  \subfloat[][Cup-like geometry] 
  { \label{fig:errorCup}
     \includegraphics[width=0.45\textwidth]{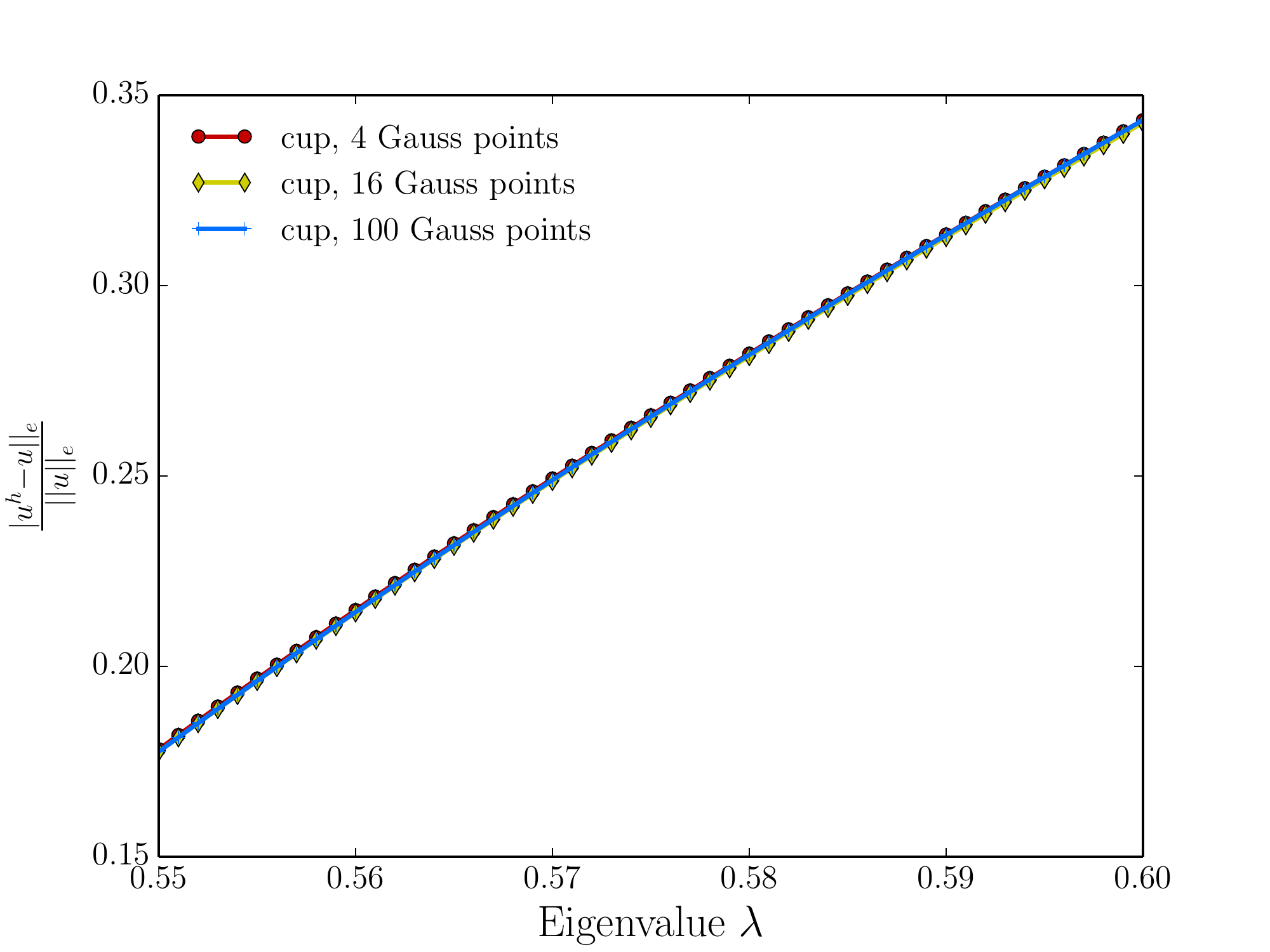}
  }\\
 \subfloat[][Saddle-like geometry] 
  { \label{fig:errorSaddle}
     \includegraphics[width=0.45\textwidth]{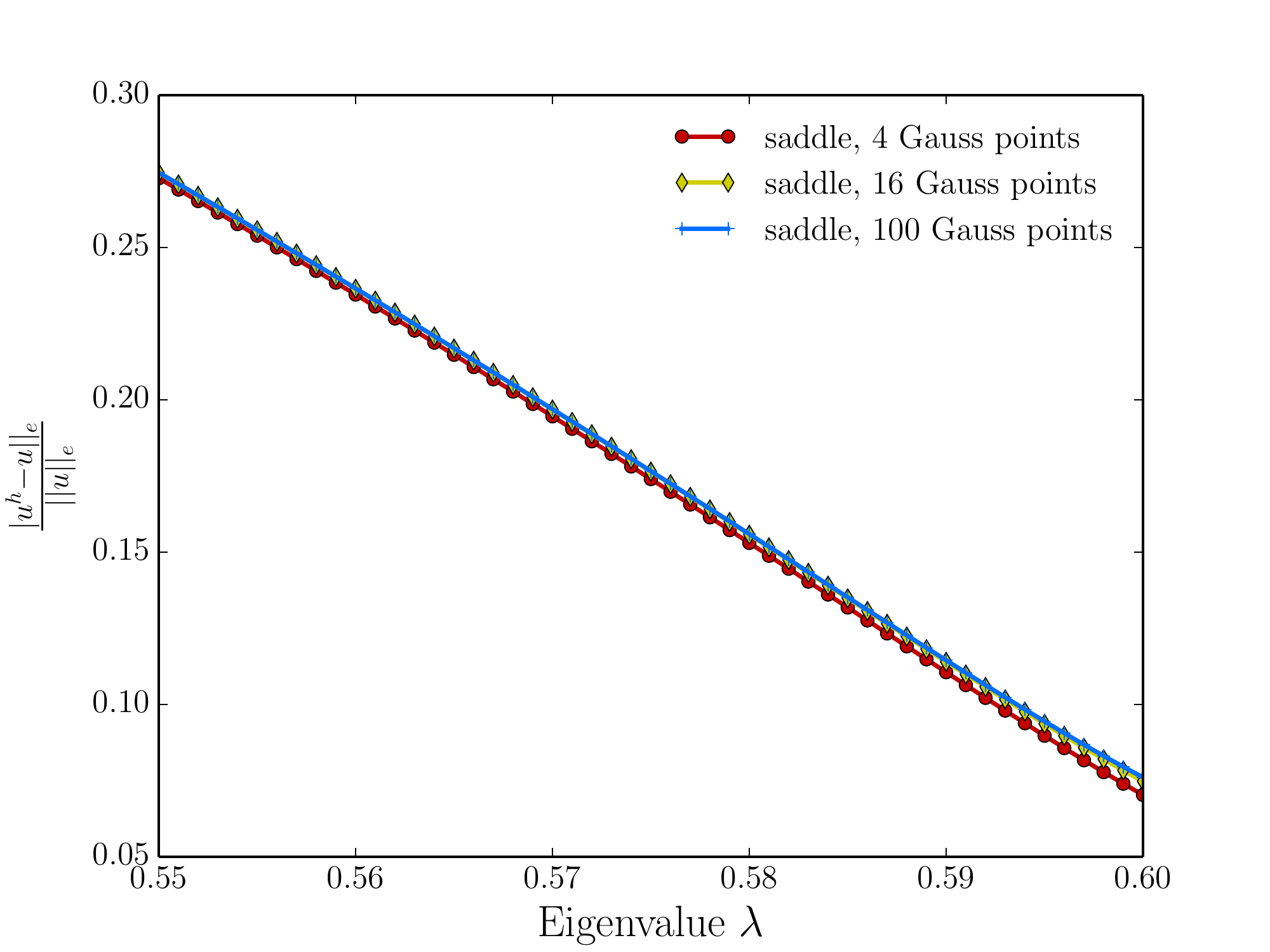}
  }
     \caption{Relative energy norm error in dependence of  the sub-dominant eigenvalue $\lambda$ and number of integration points. }
\end{figure}

\begin{figure}[ht!]
     \centering
     \includegraphics[width=0.45\textwidth]{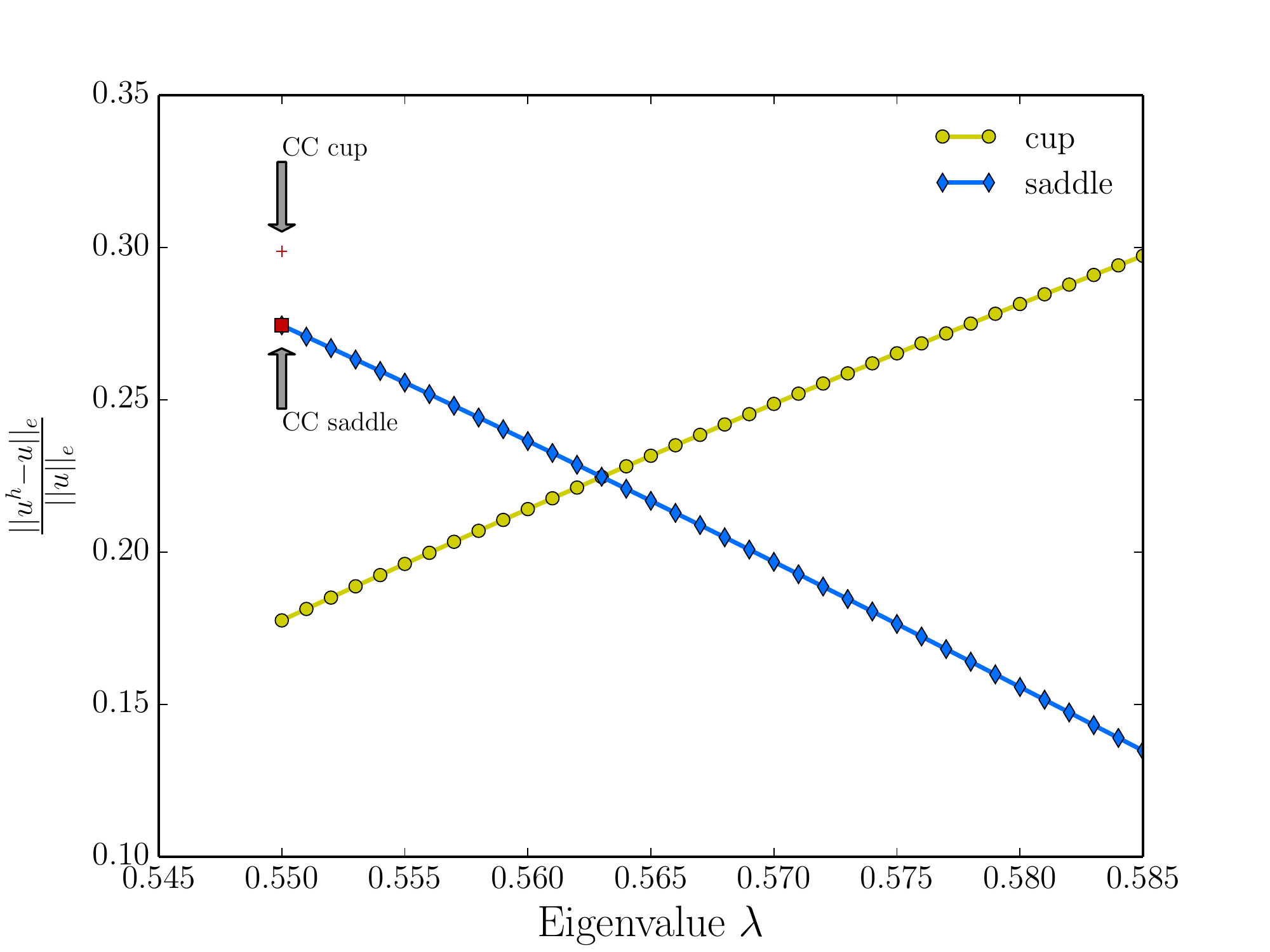}
    \caption{Relative energy norm error  for cup- and saddle-like geometries in dependence of the sub-dominant eigenvalue $\lambda$ and for the Catmull-Clark(CC) scheme. The number of integration points used for all data values is $4 \times 4$. The most optimal value for the cup-like geometry is $\lambda=0.550$, while for the saddle-like geometry it is $\lambda=0.585$. See Table~\ref{tab:weights} for the values of the optimised weights.  \label{fig:errorCupAndSaddle}}
\end{figure}

Figures~\ref{fig:errorCup} and~\ref{fig:errorSaddle} show the relative energy norm errors in approximating cup- and saddle-like geometries, respectively, when the subdominant eigenvalue $\lambda$ and number of Gauss integration points are varied. It can be seen that while $\lambda$ has a significant influence on the error the number of integration points appears to be irrelevant. Figure~\ref{fig:errorCupAndSaddle} shows the relative energy norm error both in cup- and saddle-like geometries when $4 \times 4$ integration points are used.  In comparison to Catmull-Clark weights, also indicated in Figure~\ref{fig:errorCupAndSaddle}, for $\lambda \in [0.550,0.585]$ the optimised  subdivision weights lead to a reduction of errors in both cup- and saddle-like geometries. Moreover, the most optimal value for the cup-like geometry is $\lambda=0.550$ and for the saddle-like geometry is $\lambda=0.585$, see Table~\ref{tab:weights} for the values of the optimised weights. According to Peters and Reif \cite{peters2004shape}, the obtained subdivision surfaces, same as original Catmull-Clark scheme, are $C^1$-continuous at extraordinary vertices and $C^2$ everywhere else, because the eigenvalues satisfy the required relations and the characteristic map is regular and injective.
\begin{table}[!ht]
\centering
\caption{Optimised weights in Catmull-Clark subdivision scheme for valence $v=5$ vertices.}
\label{tab:weights}
\begin{tabular}{ |c | c c c  | c | }
\hline
  & $\alpha$ & $\beta$ & $\gamma$ & $\lambda_1 = \lambda_2$ \\
\hline
Cup & 13.4575 &  0.999938 & 0.999938 & 0.550  \\
Saddle & 13.9851 & 0.824885 & 0.824885 & 0.585 \\
Original~\cite{Catmull:1978aa} & 15 & 1 & 1 & 0.550 \\
\hline
\end{tabular}		
\end{table}	

%

\subsection{Application-dependent choice of refinement weights}
\label{sec:chooseParameters}
%
When subdivision surfaces are used for finite element analysis, the solution field has quite often a mixture of both cup- and saddle-like components. And the solution field at a specific extraordinary vertex is only known after the finite element analysis.  Therefore, in a first step we use the optimal weights for the cup-like  geometry to obtain an initial finite element solution. Afterwards, for each extraordinary vertex with valence $v\geq5$,  a local shape decomposition is performed to determine whether the local solution is cup or saddle dominated. If the cup component dominates, the optimal weights for the cup-like geometry are chosen.   If instead the saddle component dominates,  the optimal weights for the saddle-like geometry are chosen. After the optimal weights for each extraordinary vertex are chosen, a second finite element analysis is performed to obtain the final solution with smaller discretisation errors. 

For thin-plate and thin-shell finite element problems the decomposition of a solution  into cup- and saddle-like components may be accomplished as described in the following. Suppose that the local finite element solution in a 3-neighbourhood of an extraordinary vertex is denoted as  $\vec p^{h}$ and has the dimensions $(12v+1) \times 3$. In a coordinate system centred at the limit position of the extraordinary vertex, according to ~\eqref{eq:limit3RingT}, we can write 
\begin{equation} \label{eq:solutionSeries}
	\vec p^h = \sum_{j=1}^{12 v}\vec{r}_j \langle \vec l_j, \vec p^h \rangle  \, . 
\end{equation}
The corresponding limit surface has at the extraordinary vertex the normal vector $\vec n \in \mathbb R^3$, defined by
\begin{equation}
	\vec n = \frac{\langle \vec l_1 , \vec p^h \rangle \times \langle \vec l_2 , \vec p^h \rangle } {| \langle \vec l_1 , \vec p^h \rangle \times \langle \vec l_2 , \vec p^h \rangle |} \, ,
\end{equation}
where the vectors $\langle \vec l_1 , \vec p^h \rangle $ and $\langle \vec l_2 , \vec p^h \rangle$ represent the two, usually non-orthogonal, tangent vectors. Multiplying~\eqref{eq:solutionSeries} with the normal vector gives by eliminating its first two terms  
\begin{equation} 
	\vec p^h \vec n^\trans= \sum_{j=3}^{12 v}\vec r_j \langle \vec l_j, \vec p^h \rangle \vec n^\trans \, . 
\end{equation}
 The  vector $\vec p^h \vec n^\trans$ represents the out-of-plane coordinates of the control points in a coordinate system aligned with the tangent plane at the extraordinary vertex. The corresponding limit surface over the characteristic domain has the following representation: 
\begin{equation}
	\label{eq:shapeDecompose1}
	\begin{split}
    u^h (\xi_1, \xi_2) = & ~\vec N^\trans \left ( \chi^{-1}(\xi_1, \xi_2)  \right ) \vec p^h \vec n^\trans \\
    = & ~\vec N^\trans \left ( \chi^{-1} (\xi_1, \xi_2) \right ) \left ( \vec r_3 \langle \vec l_3,\vec p^h \vec n^\trans \rangle  \right. \\
     &    \left. +   ~\vec r_4 \langle \vec l_4,\vec p^h \vec n^\trans \rangle   + \vec r_5 \langle \vec l_5,\vec p^h \vec n^\trans \rangle  + \cdots  \right ) 
	\end{split}
\end{equation}
and can be approximated with quadratic functions $u_3(\xi_1, \xi_2) = \xi_1^2 + \xi_2^2 \,$, $u_4(\xi_1, \xi_2) = \xi_1^2 - \xi_2^2$ and $u_5(\xi_1, \xi_2)= 2 \xi_1 \xi_2\,$, see~\eqref{eq:wSeries},  such that
\begin{equation}
	\label{eq:shapeDecompose2}
	\begin{split}
    u^h (\xi_1, \xi_2) 
     \approx &  ~\left ( \frac{\lambda_3}{\lambda_1^{2}} \right ) ^\ell \frac{\langle \vec l_3,\vec p^h \vec n^\trans \rangle}{\langle \vec l_3, u_3(\vec r_1, \vec r_2)\rangle}(\xi_1^2+\xi_2^2)  \\
    & + \left ( \frac{\lambda_4}{\lambda_1^{2}} \right )^\ell  \frac{\langle \vec l_4,\vec p^h \vec n^\trans \rangle}{\langle \vec l_4, u_4(\vec r_1, \vec r_2)\rangle}(\xi_1^2 - \xi_2^2)   \\
    & + \left ( \frac{\lambda_5}{\lambda_1^{2}} \right )^\ell \frac{\langle \vec l_5,\vec p^h \vec n^\trans \rangle}{\langle \vec l_5, u_5(\vec r_1, \vec r_2)\rangle}(2 \xi_1 \xi_2) + \cdots \\   
     \coloneqq &  ~k_3 (\xi_1^2+\xi_2^2) + k_{4} (\xi_1^2-\xi_2^2) + k_{5} (2\xi_1\xi_2) + \cdots  \, , 
	\end{split}
\end{equation}
where $\ell$ denotes the refinement level required to evaluate at the point $(\xi_1, \xi_2)$ using the Stam~\cite{Stam:1998aa} algorithm. The refinement level dependent factors always vanish when, as required for curvature continuity, $\lambda_j=\lambda_1^2 $ for $j \in {3, 4, 5}$. After computing the energy densities, i.e. the integrand in~\eqref{eq:thinPlateEnergy}, for each of the three quadratic components their ratio can be determined. Moreover, the two last terms with saddle-like geometries  are energetically equivalent so that their components can be combined. This gives the following ratio between cup- and saddle-like energies: 
\begin{equation}
	R =  \ \frac{k_3^2 \|u_3\|_e^2}{k_{4}^2 \|u_4\|_e^2+k_{5}^2 \|u_5\|_e^2} = \frac{k_3^2}{(k^2_4 + k_5^2)} \frac{1+\mu}{1-\mu} \, .
\end{equation}
In numerical computations the ratio $R$ is used to decide which set of subdivision refinement weights to use.

%
\section{Examples \label{sec: examples}}
%
We consider the finite element analysis of thin plates to demonstrate the benefits of the optimised subdivision weights over  Catmull-Clark  weights. The plates are square shaped, simply supported and subjected to either uniform or sinusoidal distributed transversal loads, see Table~\ref{tab:plates}. The thin plate energy functional depends on the second derivatives of the displacement field, c.f.~\eqref{eq:thinPlateEnergy}.  Hence, the accurate approximation of the quadratic terms in the solution field is crucial. For details of finite element implementation we refer to Cirak et al.~\cite{Cirak:2000aa, Cirak:2011aa}.  The analytical solutions of all the computed problems are known and can be found in Timoshenko et al.~\cite[Chapter~5]{Timoshenko1959}. 
\begin{table}[!ht]
	\centering
	\caption{Geometry, material and loading of the computed thin plates. \label{tab:plates}}	
	\begin{tabular}{ |  m{3.5cm}   l |   }
		\hline
 		Length & $L_x = 10$, $L_y=10$ \\ 
 		Thickness & $t=0.1$ \\  
 		Young's modulus & $E =200 \times 10^9 $ \\
 		Poisson's ratio & $\mu = 0.3$ \\
 		Uniform loading & $p =  10^4 $ \\
 		Sinusoidal loading & $p_s = p \sin(2\pi x/L_x)\sin(2\pi y/L_y) $ \\
 	\hline     
\end{tabular}

\end{table}

\begin{figure}[ht!]
  \centering
  \subfloat[][Symmetric mesh] 
  {
  	\includegraphics[width=0.21\textwidth]{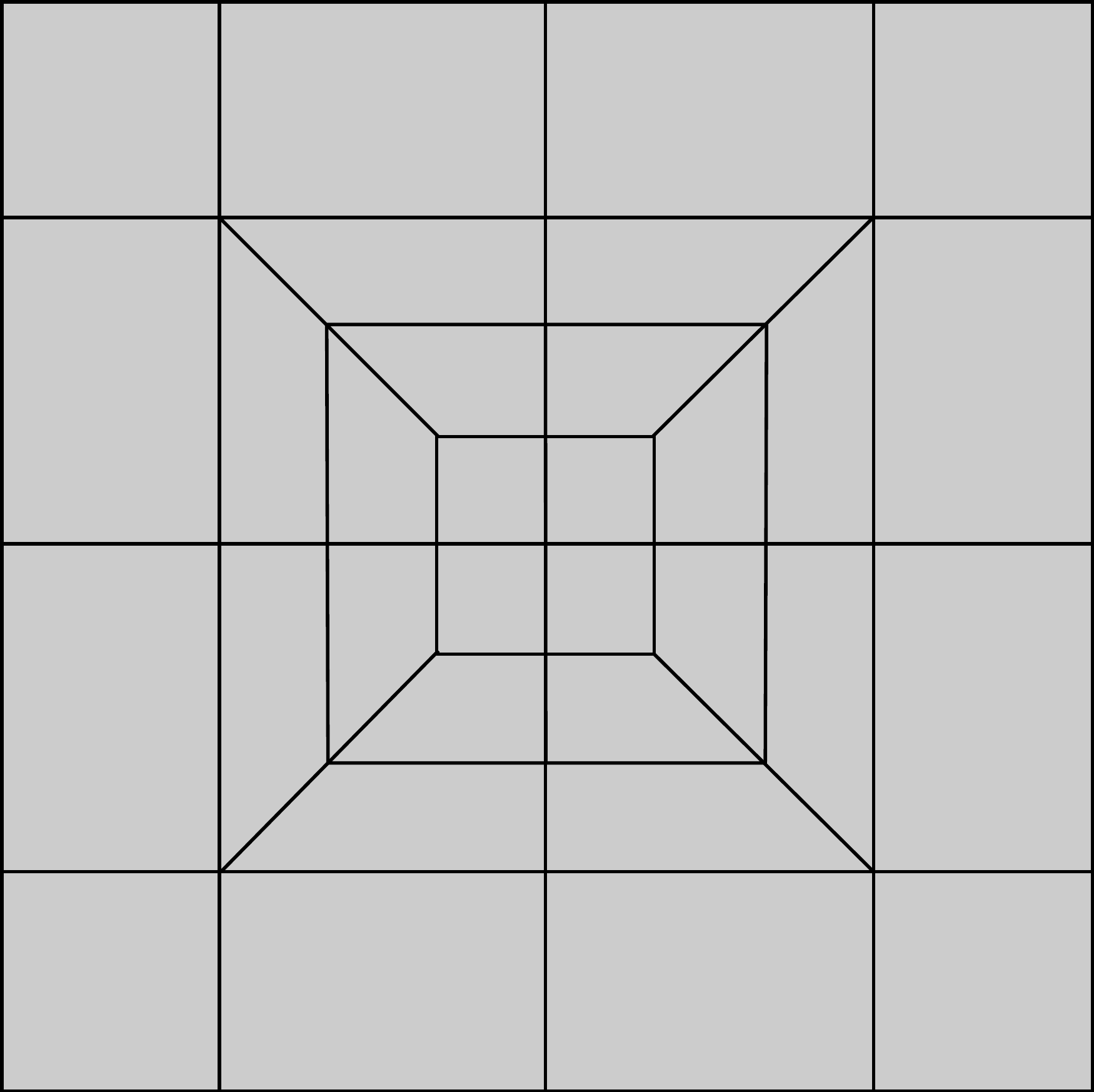}
  	\label{fig:symmetricMesh}
  }
  \hfil
  \subfloat[][Asymmetric mesh] 
  {
  	\includegraphics[width=0.21\textwidth]{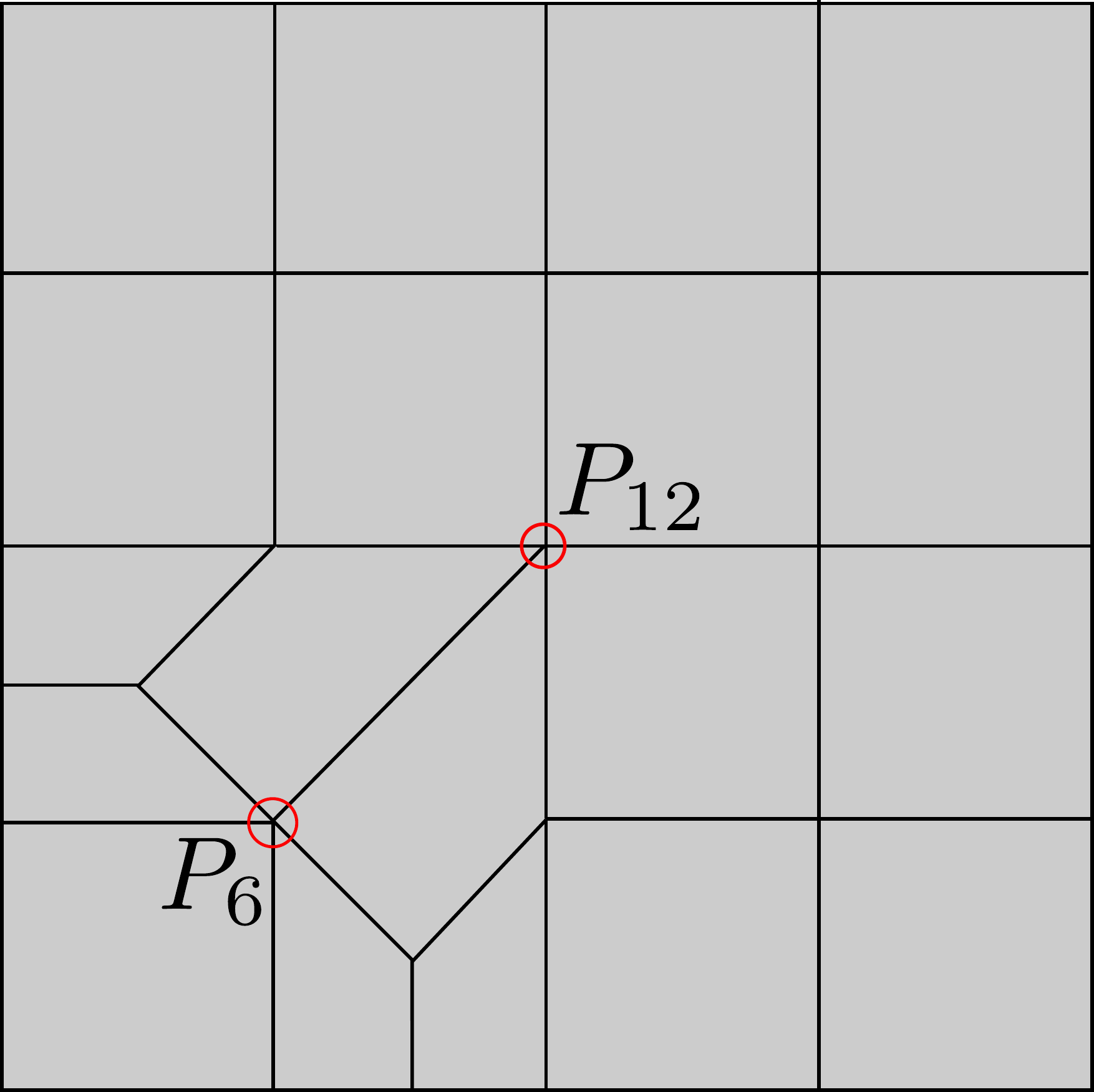}
  	\label{fig:asymmetricMesh}
  } 
  \caption{Initial unstructured coarse control meshes with level $\ell =0$. }
  \label{fig:controlMeshes}
\end{figure}

Two different unstructured control meshes shown in Figure~\ref{fig:controlMeshes} are used in the numerical computations. Both meshes have extraordinary vertices with valences $v=3$ and $v=5$. We do not use optimised subdivision weights for valence $v=3$ vertices with subdominant eigenvalues $\lambda_1 = \lambda_2 <  0.5$. During subdivision refinement their 1-neighbourhoods shrink faster than the  1-neighbourhoods of other vertices with $\lambda_1 = \lambda_2 \ge 0.5$, see also Figure~\ref{fig:characteristicMesh}. Hence, it can be expected that the benefits of optimising the weights of vertices with $v=3$ will be negligible.  For each valence $v=5$ vertex  the optimal subdivision weights are chosen independently based on the dominant component of the quadratic shape at the vertex.  

We demonstrate the benefits  of the optimised subdivision weights  over the Catmull-Clark weights by plotting the convergence of  $L_2$ and  energy norm errors. The successively refined  meshes are obtained by subdivision using the  Catmull-Clark  weights. In all examples $4 \times 4$ Gauss quadrature points are used to evaluate the finite element integrals, which appears to be sufficiently accurate as  shown in Figure~\ref{fig:errorCup} and Figure~\ref{fig:errorSaddle}.  See also \cite{juttler2016} for a systematic study on numerical integration of subdivision surfaces.


\subsection{Uniform loading, symmetric unstructured mesh}
%
As the first example, we compute the deformation of a simply supported square plate subjected to uniform transversal loading. The plate is discretised with the symmetric unstructured mesh shown  in Figure~\ref{fig:symmetricMesh}. 
Since  the displacement field is usually not known prior to a finite element analysis, we use the optimal weights for cup-like shapes to solve the plate bending problem on a level $\ell=2$ control mesh. Afterwards, for  extraordinary vertices with valence $v=5$ we decompose the local displacement field energetically  to determine whether it is cup- or saddle-dominated and choose the optimal subdivision weights accordingly. For the considered uniform loading the shape decomposition shows that saddle dominates at all valence $v=5$ vertices with cup-saddle ratio $R = 0.766$. Therefore, we choose the optimal  weights for saddle-dominated shapes for all $v=5$ vertices and study the convergence of the finite element solution using meshes  from levels $\ell=1$ to  $\ell=5$. 

\begin{figure}[!ht]
  \centering
  \subfloat[][Level $\ell = 2$ control mesh] 
  {
  	\includegraphics[width=0.25\textwidth]{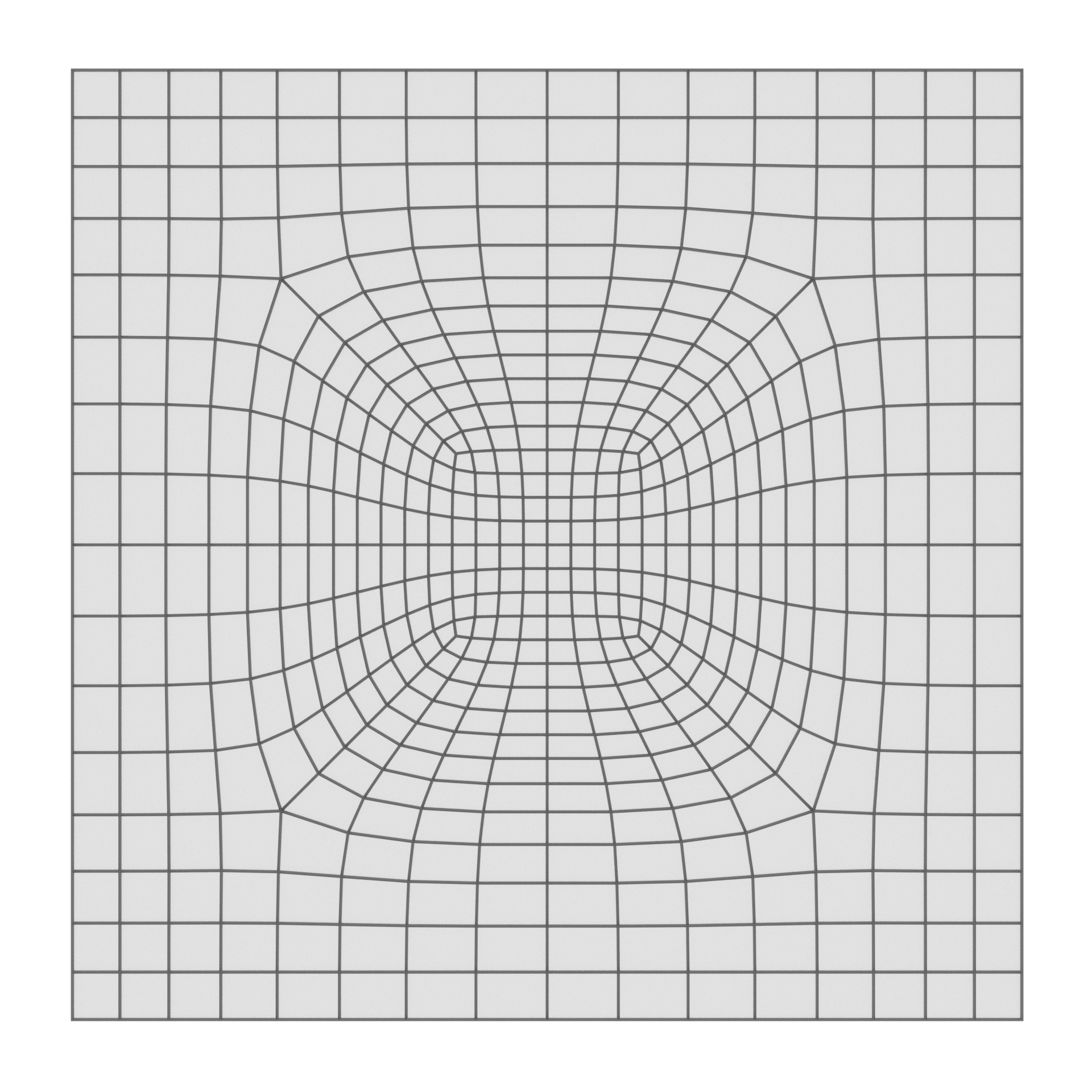}
	\label{fig:symmetricMeshL2}
  }
  \hfil
  \subfloat[][Deformed plate] 
  {
        \includegraphics[width=0.375\textwidth]{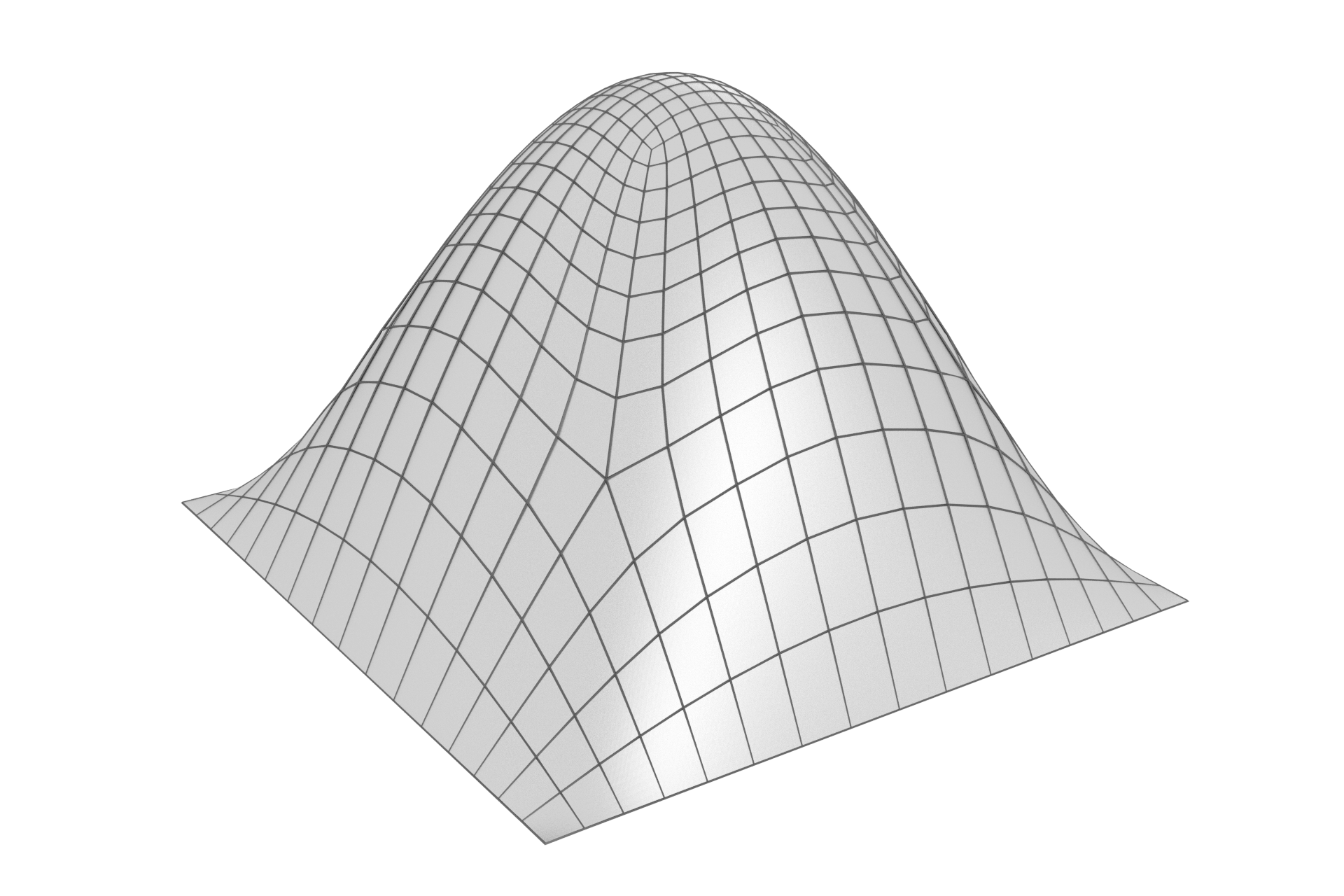}
  }
  \caption{Control mesh with fourfold symmetry and the deformed plate under uniform loading. The shown control mesh is obtained by subdividing the symmetric coarse control mesh in Figure~\ref{fig:symmetricMesh} twice using Catmull-Clark weights.}
  \label{fig:plate1-dispContour}
\end{figure}

\begin{figure}[!ht]
  \centering
  \subfloat[][$L_2$ norm error] 
  {
        \includegraphics[width=0.45\textwidth]{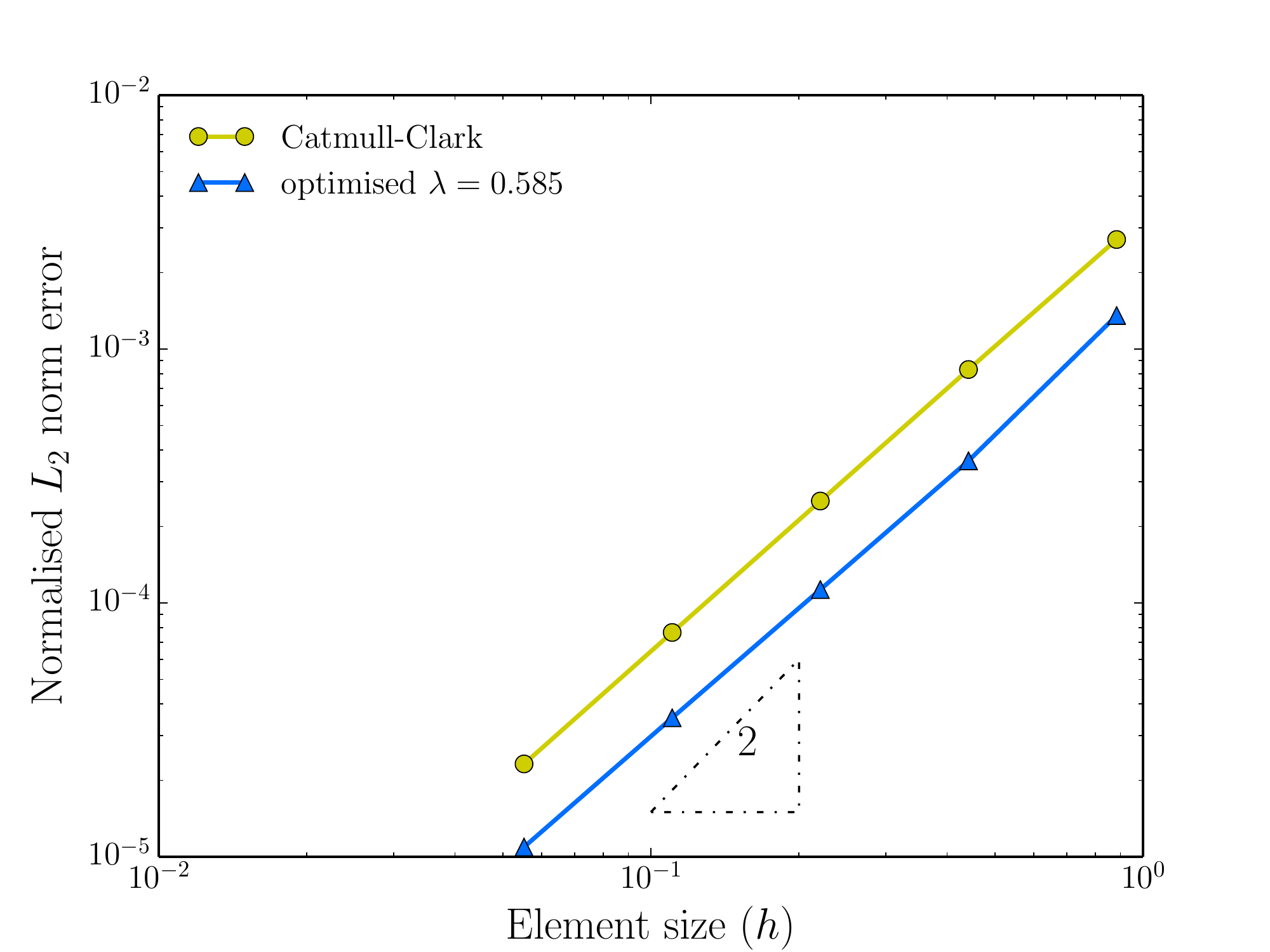}
  }\\
  \subfloat[][Energy norm error] 
  {
        \includegraphics[width=0.45\textwidth]{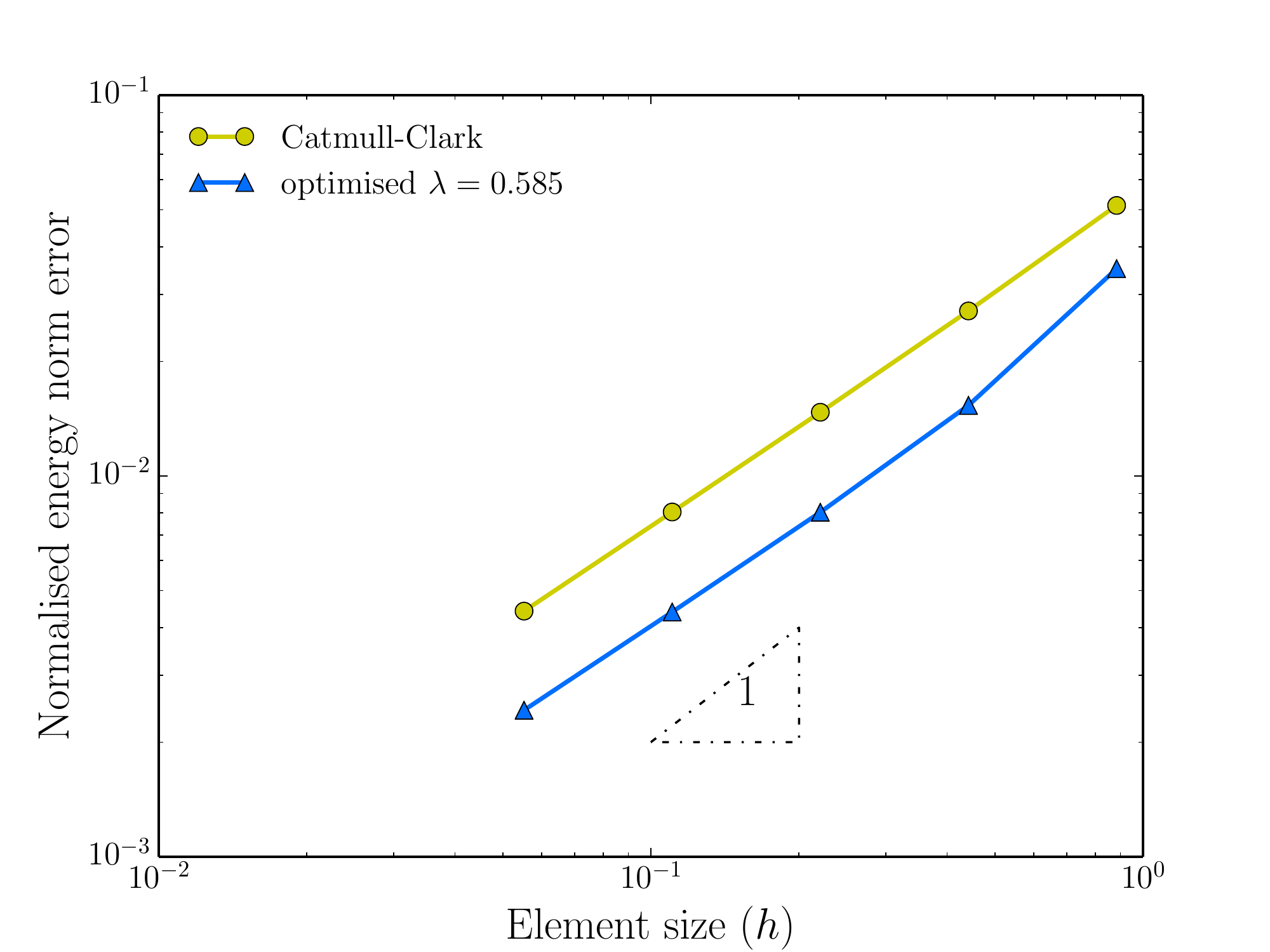}
  }
  \caption{Uniform loading with symmetric unstructured mesh. Saddle dominates ($R = 0.766$) at all valence $v=5$ vertices. 
  Optimisation reduces the~$L_2$ norm error by more than~$50 \%$ and energy norm error by more than~$45\%$. See Table~\ref{tab:weights} for the values of the optimised weights corresponding to $\lambda=0.585$.}
  \label{fig:plate1-error}
\end{figure}

Figure~\ref{fig:plate1-dispContour} shows the level $\ell=2$ control mesh and the deformed plate. The convergence of $L_2$ and energy norm errors are plotted in  Figure~\ref{fig:plate1-error}. The optimised refinement weights reduce the $L_2$ norm error by more than $50 \%$ and the energy norm error by more than $45\%$ in comparison to Camtull-Clark subdivision weights.
%

\subsection{Sinusoidal loading, symmetric unstructured mesh}
%
Next, we compute the deformation of a simply supported square plate discretised with the unstructured mesh shown in Figure~\ref{fig:symmetricMesh} and subjected to  sinusoidal loading. Compared with the first example, the only difference is that sinusoidal loading is applied  instead of a uniform loading. 
After the first finite element analysis with the optimal weights for the cup-like geometry the local shape decomposition shows that the cup component dominates at all valence $v=5$ vertices with cup-saddle ratio $R = 29.2$. Therefore, we choose the optimal weights for cup-dominated shapes for all $v=5$ vertices and study the convergence of the finite element  solution using meshes from levels $\ell=1$ to $\ell=5$. 

\begin{figure}[!ht]
  \centering
        \includegraphics[width=0.375\textwidth]{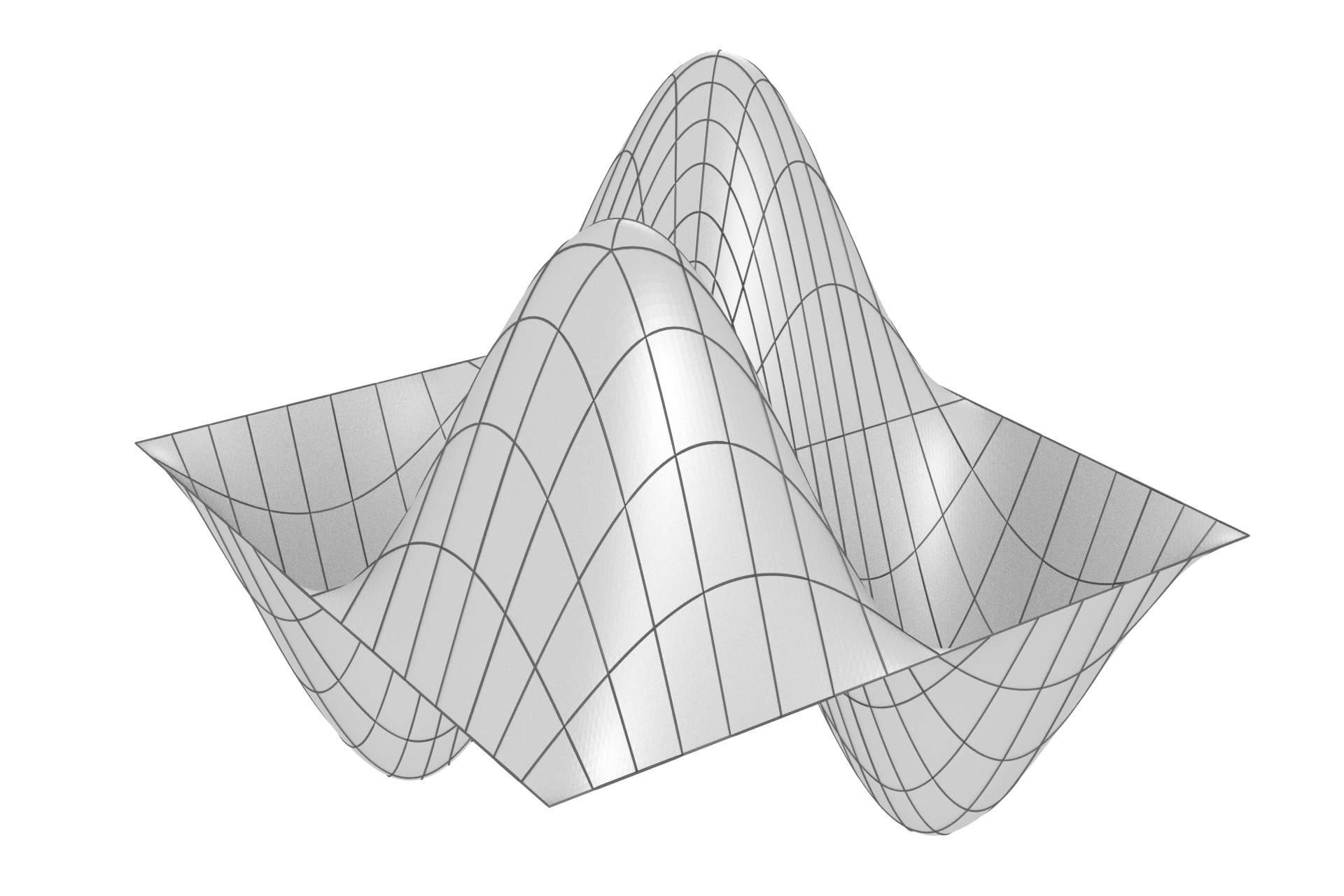}
  \caption{The deformed plate under sinusoidal loading on a level $\ell=2$ symmetric control mesh. 
  \label{fig:plate2-dispContour} }
\end{figure}

\begin{figure}[!ht]
  \centering
  \subfloat[][$L_2$ norm error] 
  {
  	\includegraphics[width=0.45\textwidth]{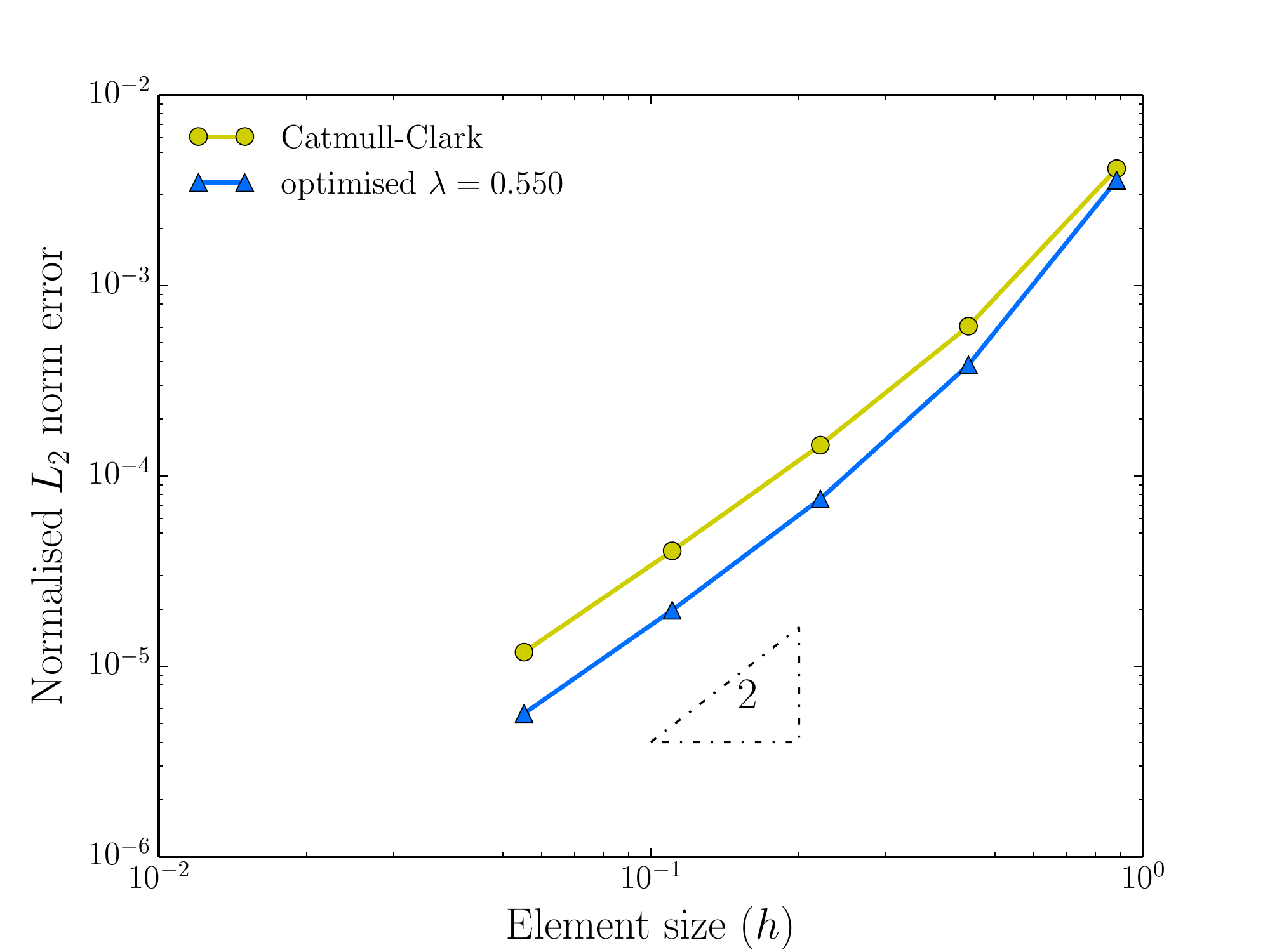}
  } \\
 \subfloat[][Energy norm error] 
  {
  	\includegraphics[width=0.45\textwidth]{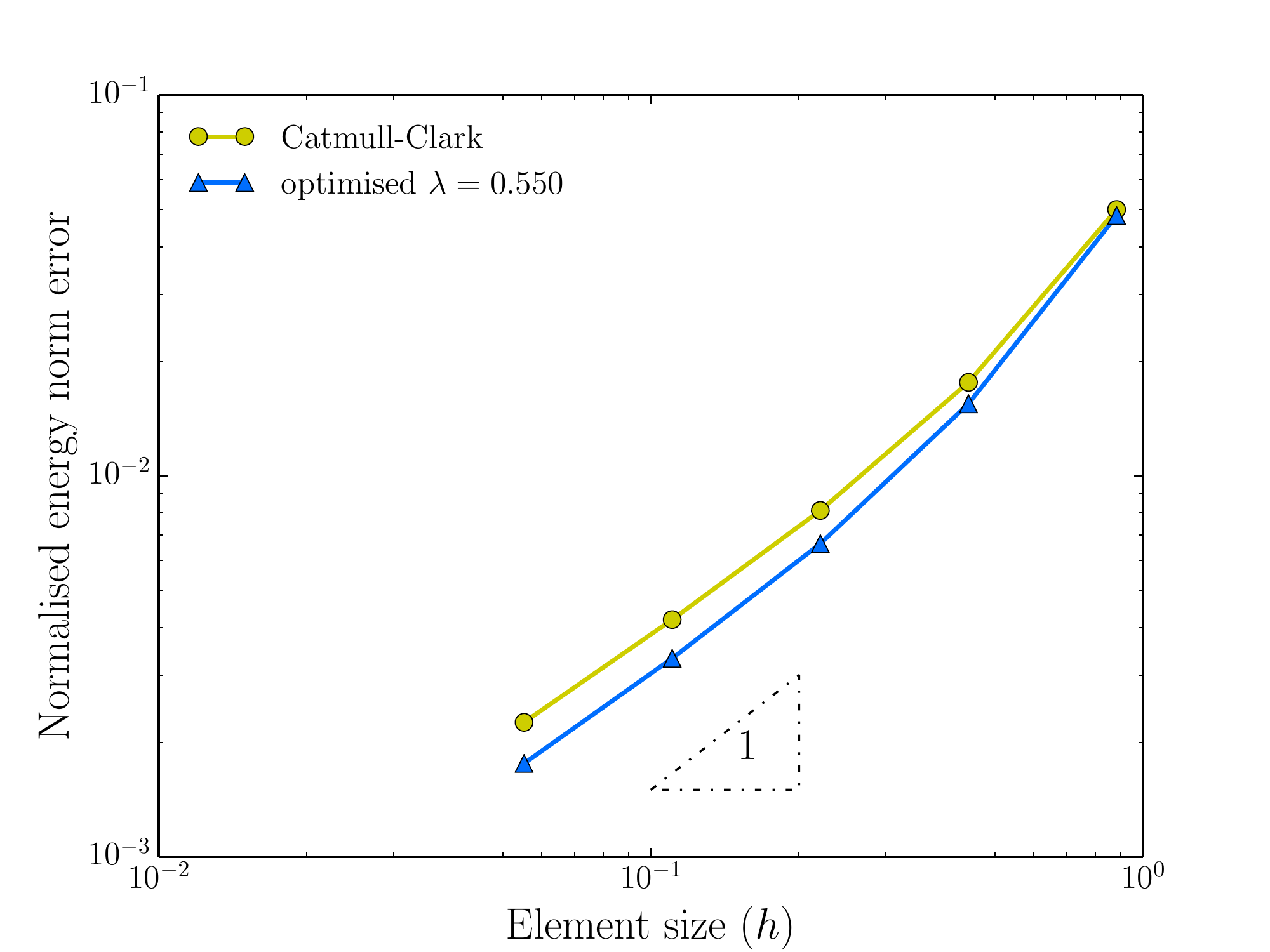}
  }
  \caption{Sinusoidal loading with symmetric unstructured mesh. Cup dominates ($R=29.2$) at all valence $v=5$ vertices.
  	Optimisation reduces the $L_2$ error by more than $50 \%$ and energy norm error by more than $20\%$. See Table~\ref{tab:weights} for the values of the optimised weights corresponding to $\lambda=0.550$.}
  \label{fig:plate2-error}
\end{figure}

Figure~\ref{fig:plate2-dispContour} shows the  level $\ell=2$ control mesh and the deformed plate. The convergence of $L_2$ and energy norm errors are plotted in  Figure~\ref{fig:plate2-error}. The optimised refinement weights reduce the $L_2$ norm error by more than $50 \%$ and the  energy norm error by more than $20\%$ in comparison to Camtull-Clark weights.

\subsection{Sinusoidal loading, asymmetric unstructured mesh} 
In this last example, we compute the deformation of a simply supported square plate discretised with the asymmetric unstructured mesh  shown in Figure~\ref{fig:asymmetricMesh} and subjected to sinusoidal loading. 
As in the first two examples, we obtain the first finite element solution using the optimal weights for cup-like shapes on a level $\ell=3$ control mesh. Subsequently, for each extraordinary vertex with valence $v=5$ the local displacement field is decomposed  to determine whether it is cup or saddle dominated  and the optimal weights are chosen accordingly. The local shape decomposition shows that the local solution at vertex $P_6$ is cup dominated with cup-saddle ratio $R=85.8$ and at vertex $P_{12}$ it is saddle dominated with $R= 0.0463$. We choose optimised cup weights for $P_{6}$ and saddle weights for $P_{12}$ and study the convergence of the finite element solution using meshes from levels $\ell=1$ to $\ell=5$.

\begin{figure}[!ht]
  \centering
  \subfloat[][Level-3 control mesh] 
  {
  	\includegraphics[width=0.25\textwidth]{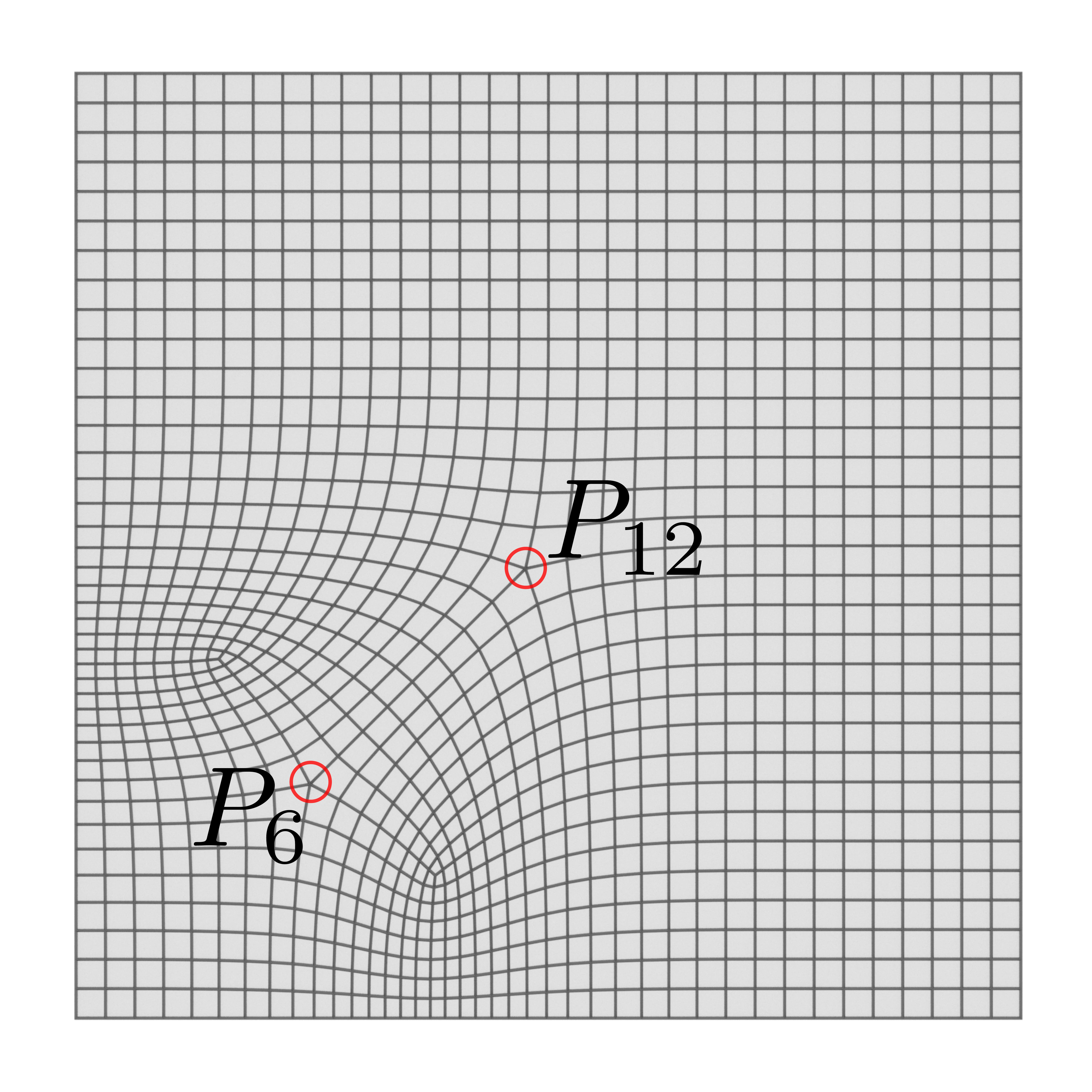}
  }
  \hfil
  \subfloat[][Deformed plate] 
  {
        \includegraphics[width=0.375\textwidth]{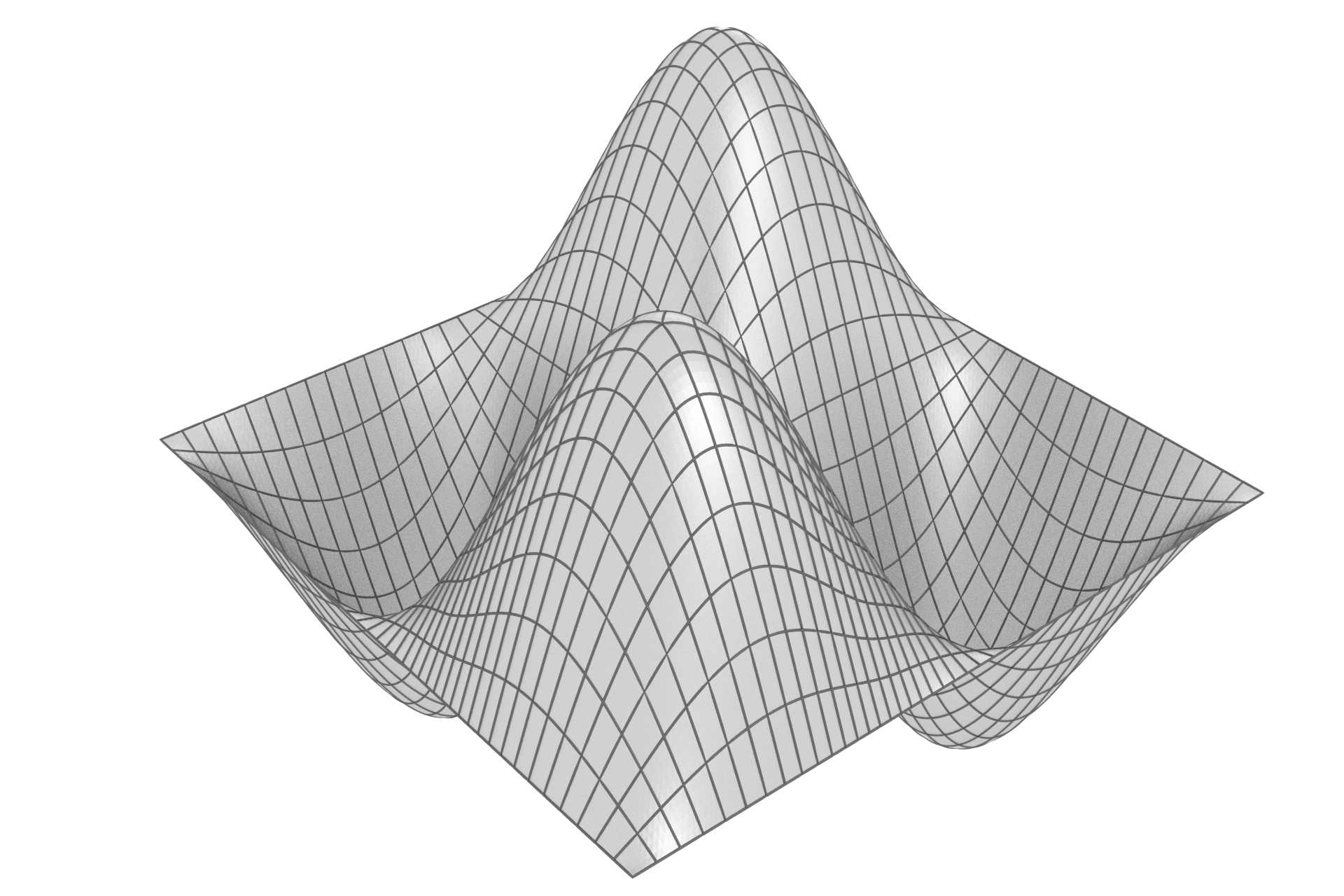}
  }
  \caption{Asymmetric control mesh and the deformed plate under sinusoidal loading. The control mesh is obtained by subdividing the asymmetric coarse control mesh shown in Figure~\ref{fig:asymmetricMesh} three times with Catmull-Clark subdivision weights.}
  \label{fig:plate3-disp}
\end{figure}

\begin{figure}[!ht]
  \centering
  \subfloat[][$L_2$ norm error] 
  {
  	\includegraphics[width=0.448\textwidth]{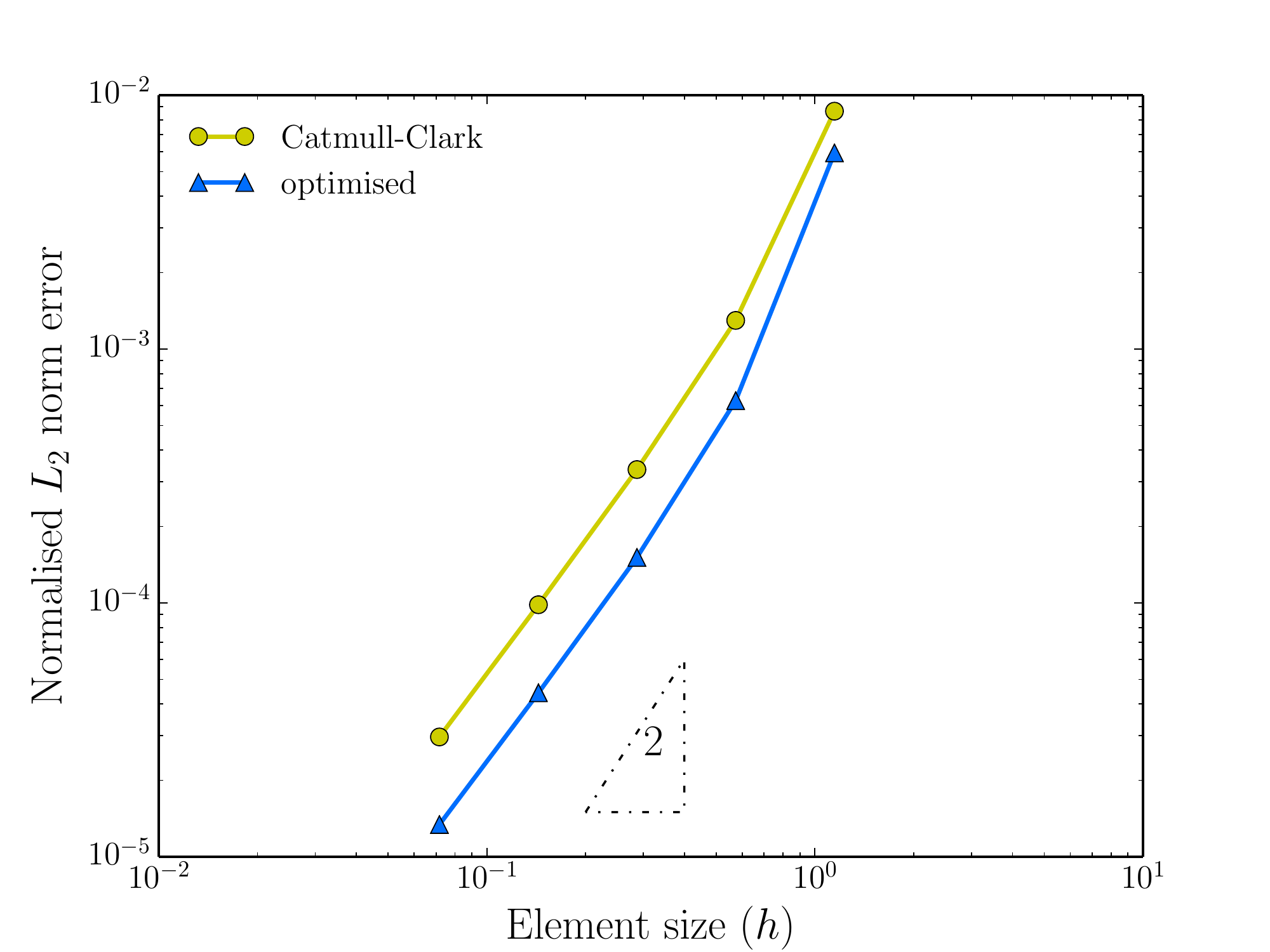}
  } \\
 \subfloat[][Energy norm error] 
  {
  	\includegraphics[width=0.448\textwidth]{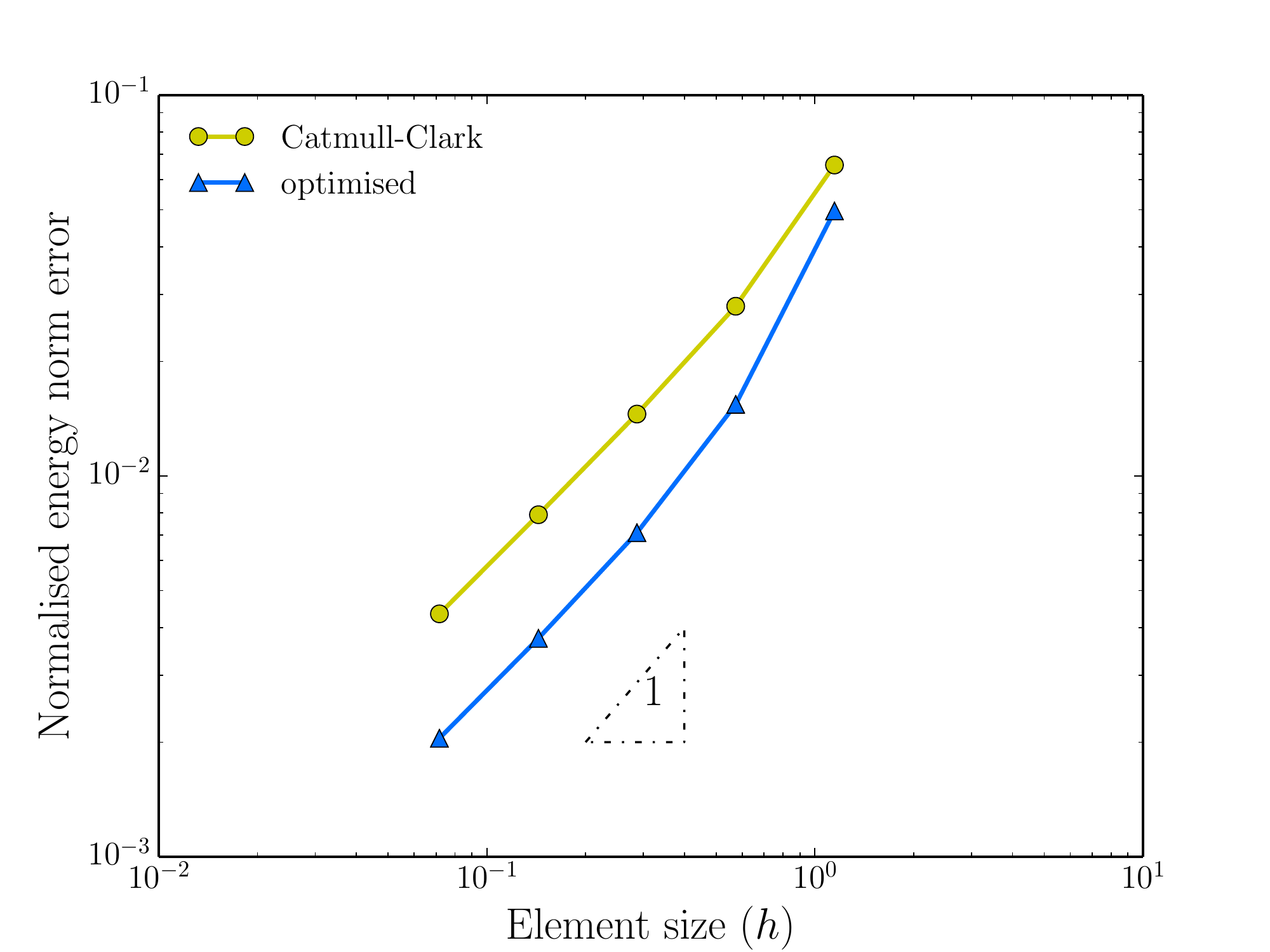}
  }
  \caption{Sinusoidal loading with asymmetric unstructured mesh. Cup dominates ($R=85.8$) at vertex $P_6$ and saddle dominates ($R=0.0463$) at vertex $P_{12}$. Therefore, $\lambda=0.550$ is chosen for vertex $P_6$ and $\lambda=0.585$ is chosen for vertex $P_{12}$. Tuning reduces both the $L_2$ error and energy norm error by more than $50 \%$. See Table~\ref{tab:weights} for the values of the optimised weights.}
  \label{fig:plate3-error}
\end{figure}

Figure~\ref{fig:plate3-disp} shows level $\ell=3$ control mesh and the deformed plate. The convergence of $L_2$ and energy norm errors are plotted in  Figure~\ref{fig:plate3-error}. The optimised refinement weights reduce both $L_2$ error and energy norm error by more than $50\%$ in comparison to Catmull-Clark weights.

\section{Conclusions \label{sec: conclusions}}
%
We have shown that significant reductions in discretisation errors in $L_2$ and energy norms can be achieved when subdivision weights around extraordinary vertices are optimised for finite element analysis.  Although this was demonstrated for Catmull-Clark subdivision surfaces, a similar approach can be developed for other subdivision schemes as well. During finite element analysis the subdivision weights at each extraordinary vertex are chosen depending on whether the local solution has a more cup- or saddle-like shape. Two sets of weights, one for cup and the other for saddle, were derived which depend only on the valence of the extraordinary vertex. We discussed only valence $v=5$ because it is, in addition to $v=3$, one of the most occurring valences for quad meshes. The same implementation applies to extraordinary vertices with $v > 5$ with no further modification.
For the case valence $v=3$ we observed that the improvement  is not significant. This is as to be expected given that the 1-neighbourhood of  a valence $v=3$ vertex shrinks faster than of other valences. In the optimisation process three of the subdivision weights, $\alpha$, $\beta$ and $\gamma$ in the 1-neighbourhood of an extraordinary vertex were selected as degrees of freedom. By considering the 2-neighbourhood of a vertex it would have been possible to optimise more than three subdivision weights. This may lead to even larger reductions in the errors although the considered optimisation problems become larger. Finally, subdivision surfaces are equally well suited for finite element analysis and modelling  of geometries with arbitrary topology. With the derived optimised weights, geometric models created with subdivision surfaces in the new engineering design systems can be analysed much more efficiently.

\section*{Acknowledgement}
Partial support through Trimble Inc and Cambridge Trust is gratefully acknowledged.

\bibliographystyle{elsarticle-num-names}
\bibliography{tunedSubdivision}

\end{document}